\documentclass{amsart}
\usepackage{graphicx}
\newtheorem{thm}{Theorem}[section]
\newtheorem*{thmnonnum}{Theorem}
\newtheorem*{propnonnum}{Proposition}
\newtheorem*{lemnonnum}{Lemma}

\newtheorem{cor}[thm]{Corollary}
\newtheorem{lem}[thm]{Lemma}
\newtheorem{prop}[thm]{Proposition}
\newtheorem{defi}[thm]{Definition}
\theoremstyle{remark}
\newtheorem{domanda}[thm]{Question}

\newenvironment{demo}{\noindent{\bf Proof: }}{\hfill$\Box$\medskip}
\numberwithin{equation}{section}

\newcommand\NN {{\mathbb N}}
\newcommand\ZZ {{\mathbb Z}}
\newcommand\QQ {{\mathbb Q}}
\newcommand\RR {{\mathbb R}}
\newcommand\CC {{\mathbb C}}

\newcommand\HH {{\mathbb H}}
\newcommand\PP {{\mathbb P}}

\newcommand\soduer{{\rm SO(2,\RR)}}
\newcommand\glduer{{\rm GL(2,\RR)}}
\newcommand\slduez{{\rm SL(2,\ZZ)}}
\newcommand\sldz{{\rm SL(d,\ZZ)}}
\newcommand\slduer{{\rm SL(2,\RR)}}

\newcommand\pslduez{{\rm PSL(2,\ZZ)}}
\newcommand{\diff}{\texttt{Diff}}

\newcommand\angolo{{\rm angle }}
\newcommand{\sys}{\texttt{Sys}}
\newcommand\area{{\rm Area }}
\newcommand\leb{{\rm Leb }}

\newcommand\hol{{\rm Hol }}

\newcommand\modu{{\rm Mod }}

\newcommand{\trasl}{\texttt{Transl}}

\newcommand\al{\alpha}

\newcommand\ga{\gamma}
\newcommand\de{\delta}

\newcommand\la{\lambda}
\newcommand\si{\sigma}
\newcommand\te{\theta}

\newcommand\De{\Delta}

\newcommand\Th{\Theta}

\newcommand\cA{{\mathcal{A}  }}

\newcommand\cC{{\mathcal{C}  }}

\newcommand\cF{{\mathcal{F}  }}
\newcommand\cG{{\mathcal{G}  }}
\newcommand\cH{{\mathcal{H}  }}
\newcommand\cI{{\mathcal{I}  }}

\newcommand\cP{{\mathcal{P}  }}
\newcommand\cQ{{\mathcal{Q}  }}
\newcommand\cR{{\mathcal{R}  }}
\newcommand\cS{{\mathcal{S}  }}
\newcommand\cT{{\mathcal{T}  }}

\newcommand\cV{{\mathcal{V}  }}

\begin{document}

\title[Khinchin condition and asymptotic laws]{Khinchin type condition for translation surfaces and asymptotic laws for the Teichm\"uller flow}
\author{Luca Marchese}
\address{Section de math\'ematiques, case postale 64, 2-4 Rue du Li\`evre, 1211 Geneve, Suisse.}
\email{l.marchese@sns.it}

%\date{}%
%\dedicatory{}%
%\commby{}%

\begin{abstract}
We study a diophantine property for translation surfaces, defined in term of saddle connections and inspired by the classical theorem of Khinchin. We prove that the same dichotomy holds as in Khinchin' result, then we deduce a sharp estimation on how fast the typical Teichm\"uller geodesic wanders towards infinity in the moduli space of translation surfaces. Finally we prove some stronger result in genus one.
\end{abstract}

% ----------------------------------------------------------------

\maketitle

\tableofcontents

\section{Introduction}

The moduli space of flat tori is identified with the \emph{modular surface}, that is the quotient $\HH/\pslduez$ of the hyperbolic half plane by the action of \emph{Moebius transformations}, which is homeomorphic to a punctured sphere and it has finite area. A neighborhood of the cusp corresponds to those flat tori with a very short closed geodesic, or equivalently to points in the standard fundamental domain with big imaginary part. The geodesic flow $g_{t}$ acts ergodically on the unitary tangent bundle of the modular surface, therefore a generic geodetic makes infinitely many excursions to the cusp. The rate of this phenomenon is quantified by the so-called \emph{logarithmic law}.

\begin{thmnonnum}
For any point $z$ in $\HH/\pslduez$ and almost any unitary tangent vector $v$ at $z$, if $g_{t}$ is the geodetic at $z$ in the direction of $v$ we have
$$
\limsup_{t\to\infty}\frac{d(g_{t}(z),z_{0})}{\log t}=1/2,
$$
where $z_{0}$ is any point in $\HH/\pslduez$ chosen as \emph{center} and $d$ denotes the \emph{Poincar\'e distance}.
\end{thmnonnum}

The logarithmic law has been generalized in many settings, in particular Sullivan proved it for the geodesic flow on manifolds with negative curvature (see \cite{sullivan}) and Masur proved a logarithmic law for \emph{Teichm\"uller geodetics} on the moduli space of complex curves of any genus (see \cite{masurtre}).

The geodesic flow on $\HH/\pslduez$ has a well know relation with the \emph{continued fraction algorithm}, which have been described in \cite{series}, and in general with arithmetics. In particular the logarithmic law is strictly related to the Khinchin theorem (see \cite{kin}), which concerns the general diophantine condition on a real number $\al$ in $(0,1)$ defined by the equation:
\begin{equation}\label{eqkhinchine}
\{n\al\}< \varphi (n),
\end{equation}
where $\{\cdot\}$ denotes the fractionary part and $\varphi :\mathbb{N}\to \mathbb{R}_{+}$ is a positive sequence such that $n\varphi (n)$ is monotone decreasing.

\begin{thmnonnum}[Khinchin]
Let $\varphi :\mathbb{N}\to \mathbb{R}_{+}$ be a positive sequence as above.
\begin{itemize}
\item
If $\sum_{n\in \mathbb{N}} \varphi(n) <\infty $ then equation (\ref{eqkhinchine}) has just finitely many solutions $n\in \mathbb{N}$ for almost any $\al$.
\item
If $\sum_{n\in \mathbb{N}} \varphi(n)=\infty $ then for almost any $\al$ equation (\ref{eqkhinchine}) has infinitely many solutions $n\in\NN$.
\end{itemize}
\end{thmnonnum}

\subsubsection*{Translation surfaces}

The natural generalization to higher genus of a flat torus is the notion of \emph{translation surface}, that is a compact, orientable and boundary-less flat surface $X$, with conical singularities whose angle is a multiple of $2\pi$. If $g$ is the genus of the surface, $\Sigma=\{p_{1},..,p_{r}\}$ is the set of conical singularities and $k_{1},..,k_{r}$ are positive integers such that for any $i=1,..,r$ the angle at $p_{i}$ is $2k_{i}\pi$, then we have the relation $k_{1}+..+k_{r}=2g+r-1$. Flat neighborhoods in $X$ are naturally identified with open sets in $\CC$, that is they admit a local coordinate $z$, thus $X\setminus \Sigma$ inherits the structure of Riemann surface and it is easy to see that the structure extends to $X$. Since the angles at conical singularities are multiples of $2\pi$ then the holonomy group is trivial and $dz$ defines a complex one form on $X\setminus \Sigma$, which extends by $z^{k_{i}-1}dz$ at any point $p_{i}\in\Sigma$, that is it has a zero of order $k_{i}-1$. The datum of a translation surface is therefore equivalent to the datum of a Riemann surface together with a complex one form. The \emph{moduli space} of translation surfaces of fixed genus $g$ admits a \emph{stratification}, where a \emph{stratum} is the set $\cH (k_{1},..,k_{r})$ of those translation surfaces with $r$ conical singularities with fixed values $2k_{1}\pi,..,2k_{r}\pi$ for the angles. We assume that the singularities are labeled. Any $\cH(k_{1},..,k_{r})$ is a complex orbifold with $\dim_{\CC}=2g+r-1$, that is it is locally homeomorphic to the quotient of $\CC^{2g+r-1}$ by the action of a finite group, orbifold points occurs at those translation surfaces whose underling complex structure admits non-trivial automorphisms. In general strata are neither compact nor connected, their connected components have been classified in \cite{konzo}.

A \emph{saddle connection} for a translation surface $X$ is a geodesic path $\ga:(0,T)\to X$ for the flat metric such that $\ga^{-1}(\Sigma)=\{0,T\}$, that is $\ga$ starts and ends in $\Sigma$ and it does not contain any other conical singularity in its interior. If $\ga$ is a saddle connection for $X$ we define a complex number $\hol(\ga):=\int_{\ga}w_{X}$, which is called the \emph{holonomy} of $\ga$, where $w_{X}$ is the holomorphic one form associated to $X$. We call $\hol(X)$ the set of complex numbers $v=\hol(\ga)$, where $\ga$ varies among the saddle connections of $X$. It is possible to see that $\hol(X)$ contains pure imaginary elements only for translation surfaces $X$ lying in a codimension one subset of $\cH(k_{1},..,k_{r})$, thus in particular a zero measure subset. Nevertheless $\hol(X)$ is always dense in $\PP\RR^{2}$ and in particular it accumulates to the imaginary axis. We define a diophantine condition comparing the deviation from the imaginary axis of elements of $\hol(X)$ with their norm. We introduce the class of \emph{strongly decreasing functions}, that is those functions $\varphi:[0,+\infty)\to (0,+\infty)$ such that $t\varphi(t)$ is decreasing monotone. For such $\varphi$ and for an element $v$ in $\hol(X)$ we consider the condition
\begin{equation}\label{eqphiconnessionesuptrasl}
|\Re(v)|<\varphi(|v|).
\end{equation}

For any translation surface $X$ the associated one form $w_{X}$ induces a pair of parallel vector fields $\partial_{x}$ and $\partial_{y}$ on $X\setminus \Sigma$ defined by $w_{X}(\partial_{x})=1$ and $w_{X}(\partial_{y})=\sqrt{-1}$ and called respectively \emph{horizontal} and \emph{vertical} vector field. They are not complete, since their trajectories stop if they arrive at a point of $\Sigma$. In particular any point $p_{i}\in\Sigma$ is the starting point of exactly $k_{i}$ half-trajectories of $\partial_{x}$, which are called \emph{horizontal separatrices}. A \emph{frame} for a translation surface $X$ is the datum of $r$ different horizontal separatrices $(S_{1},..,S_{r})$, such that $S_{i}$ starts at $p_{i}$ for any $i\in\{1,..,r\}$. Any $X$ admits $\prod_{i=1}^{r}k_{i}$ different choices of a frame. We denote $\widehat{X}$ the datum $(X,S_{1},..,S_{r})$ of a translation surface with frame and we call $\widehat{\cH}(k_{1},..,k_{r})$ the stratum of the moduli space of translation surfaces with frame, which is a covering space of $\cH(k_{1},..,k_{r})$ with natural projection $\widehat{X}\mapsto X$. It is possible to show that the covering is non-trivial, that is $\widehat{\cH}(k_{1},..,k_{r})$ is not homeomorphic to the disjoint union of $\prod_{i=1}^{r}k_{i}$ different copies of $\cH(k_{1},..,k_{r})$ and this is equivalent to say that a continuous choice of a frame is not possible on $\cH(k_{1},..,k_{r})$ but just on $\widehat{\cH}(k_{1},..,k_{r})$.
Conceptually the construction is the same as that of the orientable double covering of a non-orientable manifold.

Let $\widehat{X}$ be a translation surface with frame whose frame is $(S_{1},..,S_{r})$. Let $p_{j}$ and $p_{i}$ be any two points in $\Sigma$ (possibly the same) and let $(m,l)$ be a pair of integers with $1\leq m\leq k_{j}$ and $1\leq l\leq k_{i}$. We define the \emph{bundle} $\cC^{(p_{j},p_{i},m,l)}(\widehat{X})$ as the set of those saddle connections $\ga$ for $X$ which start at $p_{j}$, and in $p_{i}$ and satisfy
$$
2(m-1)\pi\leq\angolo(\ga,S_{j})<2m\pi
\texttt{ and }
2(l-1)\pi\leq\angolo(\ga,S_{i})<2l\pi,
$$
where $S_{j}$ and $S_{i}$ are the horizontal separatrices in the frame starting respectively at $p_{j}$ and $p_{i}$. A choice of a frame for $X$ therefore induces a decomposition of $\hol(X)$ into subsets
$$
\hol^{(p_{j},p_{i},m,l)}(\widehat{X}):=
\{\hol(\ga);\ga\in \cC^{(p_{j},p_{i},m,l)}(\widehat{X})\}.
$$
We proved the following dichotomy.

\begin{thm}\label{teoremab}
Let $\varphi:[0,+\infty)\to (0,+\infty)$ be a function such that $t\varphi(t)$ is decreasing monotone.
\begin{enumerate}
\item
If $\int_{0}^{+\infty}\varphi(t)dt<+\infty$ then $\hol(X)$ contains just finitely many solutions of equation (\ref{eqphiconnessionesuptrasl}) for almost any $X\in \cH(k_{1},..,k_{r})$.
\item
If $\int_{0}^{+\infty}\varphi(t)dt=+\infty$ then for almost any $\widehat{X}\in \widehat{\cH}(k_{1},..,k_{r})$, for any pair of points $p_{j},p_{i}$ in $\Sigma$ and for any pair of integers $(m,l)$ with $1\leq m\leq k_{j}$ and $1\leq l\leq k_{i}$ the set $\hol^{(p_{j},p_{i},m,l)}(\widehat{X})$ contains infinitely many solutions of equation (\ref{eqphiconnessionesuptrasl}).
\end{enumerate}
\end{thm}

Translation surfaces are strictly linked to \emph{interval exchange transformations} (i.e.t. in the following), a class of maps of the interval which has been largely studied, for example in \cite{veech}, \cite{masuruno}, \cite{zorich}. In particular theorem \ref{teoremab} is a consequence of a generalization of Khinchin theorem to i.e.t.s which is proved in \cite{luca2}. Anyway the convergent part of theorem \ref{teoremab} admits a stronger version, which is proved independently from i.e.t.s with an easier argument. Such stronger statement is based on a very natural construction, which consists in fixing a translation surface $X$ and rotating its vertical direction. If $\te$ is the rotation angle, we call $X_{\te}$ the rotated translation surface. The image of the application $\te\mapsto X_{\te}$ is an embedded circle in the stratum $\cH(k_{1},..,k_{r})$ of $X$, except for orbifold points, where the image has a finite number of self-intersections. Globally we get a smooth orbit foliation of $\cH(k_{1},..,k_{r})$ under the action of $\soduer$.

\begin{prop}\label{proposizioneb}
Let $\varphi$ be a strongly decreasing function with $\int_{0}^{+\infty}\varphi(t)dt<+\infty$. Then for any $X$ in $\cH(k_{1},..,k_{r})$ and almost any $\te\in[0,2\pi)$ the set $\hol(X_{\te})$ contains finitely many solutions of equation (\ref{eqphiconnessionesuptrasl}).
\end{prop}

Since the orbit foliation of $\cH(k_{1},..,k_{r})$ under the action of $\soduer$ is smooth, the convergent part of theorem \ref{teoremab} is an immediate consequence of proposition \ref{proposizioneb}.

\subsubsection*{Teichm\"uller flow}

For any translation surface $X$ and any element $G\in\slduer$ we can define a new translation surface $GX$ whose local charts are the composition of the local charts of $X$ with $G$ (the direction of local charts is assumed from the surface to $\CC$). We have therefore an action of $\slduer$ on $\cH(k_{1},..,k_{r})$, or on its covering $\widehat{\cH}(k_{1},..,k_{r})$. The \emph{Teichm\"uller flow} $\cF_{t}$ is the action of the diagonal elements of $\slduer$ with trace $2\cosh t$.

Any translation surface $X$ admits an area form $dz\wedge d\bar{z}$ on $X\setminus \Sigma$, which defines a smooth function $X\mapsto \area(X)$ on strata. We call $\cH^{(1)}(k_{1},..,k_{r})$ the smooth hyperboloid in $\cH(k_{1},..,k_{r})$ of area one translation surfaces and $\widehat{\cH}^{(1)}(k_{1},..,k_{r})$ the corresponding hyperboloid in the covering space of translation surfaces with frame. $\slduer$ preserves the hyperboloids and thus the Teichm\"uller flow too. The fundamental result of Masur and Veech (see \cite{masuruno} and \cite{veech}) says that $\cF_{t}$ preserves an unique smooth probability measure $\mu^{(1)}$ on $\cH^{(1)}(k_{1},..,k_{r})$ and its restriction to any connected component of the stratum is ergodic.

Strata are non compact and area one hyperboloids are not compact too, therefore by recurrence a generic orbit of $\cF_{t}$ makes excursions \emph{to infinity}. A criterion to describe non-compact part of strata can be stated in terms of the \emph{systole function} $X\mapsto\sys(X)$, which assigns to $X$ the length of its shortest saddle connection. The main property is that a sequence $X_{n}$ in $\cH(k_{1},..,k_{r})$ (or in the area-one hyperboloid) diverges, that is it escapes from any compact set, if and only if $\sys(X_{n})\to 0$. Anyway in general the non-compact part of strata is not a simple cusp, like for flat tori, and the vanishing of the systole does not give any information on the component of it towards the sequence $X_{n}$ is diverging to. In order to keep track of this information, for any $\widehat{X}$ in $\widehat{\cH} (k_{1},..,k_{r})$, for any pair of points $p_{j},p_{i}\in\Sigma$ and any pair of integers $(m,l)$ with  $1\leq m\leq k_{j}$ and $1\leq l\leq k_{i}$ we set
$$
\sys^{(p_{j},p_{i},m,l)}(\widehat{X}):=
\min \{|v|;v\in \hol^{(p_{j},p_{i},m,l)}(\widehat{X})\},
$$
that is we consider the length of the shortest saddle connection in the bundle $\cC^{(p_{j},p_{i},m,l)}(\widehat{X})$. By recurrence the generic orbit is neither divergent nor bounded, thus for any datum $(p_{j},p_{i},m,l)$ as before the quantity $\sys^{(p_{j},p_{i},m,l)}(\cF_{t}\widehat{X})$ stays bounded away from $0$ for the most of the time, but there are arbitrary big instants $t_{n}$ such that $\sys^{(p_{j},p_{i},m,l)}(\cF_{t_{n}}\widehat{X})$ becomes small. We proved the following theorem, which gives a sharp quantitative description of this phenomenon.

\begin{thm}\label{teoremac}
Let $\psi:[0,+\infty)\to(0,+\infty)$ be a monotone decreasing function.
\begin{itemize}
\item
If $\int_{0}^{+\infty}\psi(t)dt<+\infty$ then for almost any $X\in\cH^{(1)}(k_{1},..,k_{r})$ we have:
$$
\lim_{t\to\infty}
\frac
{\sys(\cF_{t}X)}
{\sqrt{\psi(t)}}
=\infty.
$$
\item
If $\int_{0}^{+\infty}\psi(t)dt=+\infty$ then for almost any $\widehat{X}\in \widehat{\cH}^{(1)}(k_{1},..,k_{r})$, for any pair of points $p_{j},p_{i}\in\Sigma$ and any pair of integers $(m,l)$ with $1\leq m \leq k_{j}$ and $1\leq l\leq k_{i}$ we have
$$
\liminf_{t\to\infty}
\frac
{\sys^{(p_{j},p_{i},m,l)}(\cF_{t}\widehat{X})}
{\sqrt{\psi(t)}}
=0.
$$
\end{itemize}
\end{thm}

For the one parameter family $\psi(t):=\min\{1,t^{-r}\}$ with $r\geq 1$ a straightforward application of both the divergent and the convergent parts of theorem \ref{teoremac} implies the following corollary.

\begin{cor}\label{corollarioc}
For any pair of points $p_{i},p_{j}\in\Sigma$, for any pair of integers $(m,l)$ with $1\leq m \leq k_{j}$ and $1\leq l\leq k_{i}$ and for almost any $\widehat{X}\in \widehat{\cH}^{(1)}(k_{1},..,k_{r})$ we have:
$$
\limsup_{t\to\infty}
\frac
{-\log\sys^{(p_{j},p_{i},m,l)}(\cF_{t}\widehat{X})}
{\log t}
=\frac{1}{2}.
$$
\end{cor}

Corollary \ref{corollarioc} corresponds to Masur's logarithm law (see \cite{masurtre}) restricted to strata of abelian differentials. It is a less general result, since Masur's logarithmic law holds in the most general setting of \emph{quadratic differentials} (i.e. the cotangent bundle of the moduli space of complex curves). On the other hand our result contains a finer information, since we are able to prove the estimation not only for the classical systole function $\sys(X)$, but for all the functions $\sys^{(p_{j},p_{i},m,l)}(\widehat{X})$.

\begin{domanda}\label{domandainvoluzionilineari}
Theorem \ref{teoremac} is a consequence of theorem \ref{teoremab}, which is based itself on a generalization of Khinchin theorem for i.e.t.s, namely theorem \ref{teoremaa} in paragraph \ref{approximationconnectionsiet(sss)}. I.e.t.s have been generalized with \emph{linear involutions} by Danthony and Nogueira (see \cite{danno}) and the latter have been related to quadratic differentials by Boissy and Lanneau (see \cite{bola}). We ask if it is possible to extend theorem \ref{teoremaa} to linear involutions, which would bring to an extension of theorem \ref{teoremac} to the setting of quadratic differentials.
\end{domanda}

\subsubsection*{Punctured tori}

Since the orbit foliation under $\soduer$ is smooth, the divergent
part of theorem \ref{teoremab} implies that for almost any translation surface $X$ the statement in the theorem holds for almost any rotated surface $X_{\te}$. We ask what we can prove for any orbit and not just for generic ones. We have some partial result in this direction for a flat torus $X$ with $r$ punctures. Since we consider a flat torus we not have true conical singularities, that is any puncture has an euclidian neighborhood, anyway their position determines a subset $\Sigma=\{p_{1},..,p_{r}\}$ of $X$ with a non-trivial geometric information. The moduli space of the pairs $(X,\Sigma)$ is a fibre bundle over the modular surface. Any point of $\Sigma$ is the starting point of exactly one horizontal separatrix, thus a frame is intrinsically defined and a bundle of saddle connections is simply determined by the pais of points $p_{j},p_{i}$ in $\Sigma$ where the saddle connections $\ga$ respectively start and end, therefore we denoted it simply $\cC^{(p_{j},p_{i})}(X)$ and we write $\hol^{(p_{j},p_{i})}(X)$ for the set of the corresponding holonomy vectors. We remark that even if points in $\Sigma$ have a flat neighborhood, we want saddle connections not have them as interior points. We proved the following theorem.

\begin{thm}\label{teoremad}
Let $\varphi$ be a strongly decreasing function with $\int_{0}^{\infty}\varphi(t)dt=+\infty$. Then for any flat torus $X$ with $r$ punctures $p_{1},..,p_{r}$ and for almost any $\te\in [0,2\pi)$ there are at least $2r-1$ different pairs of punctures $(p_{j},p_{i})$ such that the corresponding sets $\hol^{(p_{j},p_{i})}(X_{\te})$ contain infinitely many solutions of equation (\ref{eqphiconnessionesuptrasl}).
\end{thm}

\emph{Note:}
We observe that for an arbitrary set of punctures $\Sigma$ we cannot expect to have the same result of theorem \ref{teoremad} for any pair of points. This is evident for example if $X$ is the flat torus $\CC/(2\ZZ)^{2}$ with four punctures corresponding to the points of $\ZZ^{2}/(2\ZZ)^{2}$, indeed in this case the set $\hol^{(p_{i},p_{i})}(X)$ is empty for any $p_{i}\in\Sigma$. We observe also that we do not know if our result is optimal, that is neither we are able to find generically more than $2r-1$ different pairs $(p_{j},p_{i})$ as in the statement, nor we can find counterexamples.

\subsection*{Content of this paper}

\emph{Section 2} is a brief survey of the background theory of translation surfaces and interval exchange transformations. In paragraph \ref{translationsurfaces(ss)} we give rigorous definitions for translation surfaces and their moduli space. In paragraph \ref{ietandveechconstruction(ss)} we introduce the class of interval exchange transformations and their parameter space, which is union of euclidian spaces. Then we describe \emph{Veech's construction}, which is a procedure to get a translation surface as \emph{suspension} of an i.e.t.. We also explain that the construction provides us with a family of local charts on strata of the moduli space which cover a subset of full measure.

\emph{Section 3} is devoted to the proof of theorem \ref{teoremab}. In paragraph \ref{prooftheorembconvergentcase(ss)} we prove proposition \ref{proposizioneb} and according to the discussion above this proves also the convergent case of theorem \ref{teoremab}. The argument is a simple application of the Borel-Cantelli lemma on the orbit under $\soduer$ of any translation surface $X$. In paragraph \ref{bundlesofsaddleconnections(ss)} we give the definition of translation surfaces with frame, then for these objects we define the \emph{bundles} of their saddle connections. We show that the bundles are preserved by the Teichm\"uller flow on strata of translation surfaces with frame. In \S \ref{relationveech(sss)} we give an alternative definition of bundles of saddle connections, adapted to the Veech's construction. In lemma \ref{lemconfigurazioneconnessionisellaveech} we show that the two definitions of the bundles coincide. In paragraph \ref{combdefinedsaddleconnections(ss)} we give a combinatorial construction of families of saddle connection in a required bundle. We work with Veech's construction, our construction is based on the combinatorial notion of \emph{reduced triple} for an i.e.t., which is introduced in definition \ref{deftriplaridotta}. The main result of the paragraph is proposition \ref{propconnexsellacombinatorie}. In paragraph \ref{prooftheoremBdivergentcase(ss)} we give the proof of the divergent part of theorem \ref{teoremab}. We apply theorem \ref{teoremaa}, which is a generalization of Khinchin theorem for i.e.t.s and provides us with infinitely many reduced triples satisfying equation (\ref{eqapproxconnessionitsi}), a diophantine condition related to equation (\ref{eqphiconnessionesuptrasl}). Then we apply the construction of paragraph \ref{combdefinedsaddleconnections(ss)}.

\emph{Section 4} contains the proof of theorem \ref{teoremac}. We first observe that theorem \ref{teoremab} also holds on the hyperboloid $\widehat{\cH}^{(1)}(k_{1},..,k_{r})$ of area one translation surfaces with frame with respect to the associated smooth measure $\mu^{(1)}$. In paragraph \ref{preliminaryestimations(ss)} we prove some preliminary estimations on the dilatation of lengths of saddle connection under the Teichm\"uller flow. We just look at length and we are not interested on the information on the bundles they belong to, so we work with translation surfaces without any choice of a frame. In paragraph \ref{prooftheoremc(ss)} we complete the proof of theorem \ref{teoremac}. The proof of corollary \ref{corollarioc} is omitted since it is straightforward.

\emph{Section 5} contains the proof of theorem \ref{teoremad}. In paragraph \ref{veechconstructionflattori(ss)} we recall that a flat torus $X$ without vertical closed geodesics and with a set $\Sigma$ of $r$ punctures can be represented as a suspension of a \emph{rotational} i.e.t. $T$ with $r+1$ intervals, that is an i.e.t. with just one true singularity. Then we study how the parameters describing the associated rotational i.e.t. change when the vertical direction of $X$ is rotated. The main result in the paragraph is lemma \ref{lemsingolaritatsitororotato}. In paragraph \ref{prooftheoremD(ss)} we prove theorem \ref{teoremad} applying the construction of paragraph \ref{veechconstructionflattori(ss)} and using an arithmetic result, namely theorem \ref{teoremae}, which is a generalization of Khinchin theorem for a diophantine condition defined by equation (\ref{eqteoremae}). In paragraph \ref{arithmeticresult(ss)} we prove theorem \ref{teoremae}, in particular in \S \ref{classicalcontinuedfraction(sss)} we recall some very classical results for the \emph{continued fraction algorithm}, our presentation follows a geometrical interpretation. In \S \ref{twistedcontinuedfraction(sss)} we define a family of approximations for an irrational number called \emph{twisted approximation}. The proof of the theorem is completes in \S \ref{sufficientcondition(cap5)(sss)} and \ref{endproof(cap5)(sss)}.

\subsection*{Acknowledgements}

The results in this paper were obtained in my Ph-D thesis. I would like to thank Jean-Christohpe Yoccoz for many discussions and for his help in revising this work. I am also grateful to Stefano Marmi for many discussions and to Giovanni Forni and Pascal Hubert for many precious remarks.

\section{Background}

\subsection{Translation surfaces}\label{translationsurfaces(ss)}

Let $g$ be a positive integer and consider $r$ positive integers $k_{1},..,k_{r}$ such that $k_{1}+..+k_{r}=2g+r-2$. A \emph{translation surface} is given by the following data

\begin{enumerate}
\item
A compact, boundary-less, orientable topological surface $X$ of genus $g$ with a finite subset $\Sigma=\lbrace p_{1},...,p_{r}\rbrace$ with $r$ elements.
\item
An atlas on $X\setminus\Sigma$ such that the changes of charts are translations. The direction of charts is assumed from the surface to $\CC$, if $\rho$ and $\rho'$ are two element of the atlas then $\rho'\circ\rho^{-1}(z)=z+const$. Such charts are called \emph{translation charts}.
\item
For any $p_{i}\in \Sigma$ a neighborhood $V_{i}$ of $p_{i}$, a neighborhood $W$ of $0$ in $\CC$ and a ramified covering $(V_{i},p_{i})\to (W,0)$ of degree $k_{i}$ whose branches are translation charts of the atlas.
\end{enumerate}

We call $\trasl(k_{1},..,k_{r})$ the set of translation structures as above. Any element $X$ in $\trasl(k_{1},..,k_{r})$ provides us with the following structures

\begin{itemize}
\item
A Riemann surface structure on $X$ (not only on $X\setminus \Sigma$).
\item
An holomorphic $1$-form $w_{X}$ on $X$, which in translation charts is given by $dz$. The zeros of $w_{X}$ are exactly the points of $\Sigma$, in a chart around any $p_{i}\in\Sigma$ the form $w_{X}$ is given by $z^{k_{i}-1}dz$, that is it has a zero of order $k_{i}-1$.
\item
A flat metric $g_{X}:=|dz|^{2}$ defined on $X\setminus \Sigma$, with cone singularities of total angle $2k_{i}\pi$ at any point $p_{i}\in\Sigma$.
\item
Two vertical and horizontal vector fields $\partial_{y}$ and $\partial_{x}$ on $X\setminus \Sigma$ defined by $w_{X}(\partial_{x})=1$ and $w_{X}(\partial_{y})=\sqrt{-1}$ and denoted respectively \emph{horizontal} and \emph{vertical} vector field. These fields are not complete since their trajectories stop at points in $\Sigma$. We call \emph{horizontal separatrices} the trajectories of $\partial_{x}$ starting from a point $p_{i}\in \Sigma$. A singularity $p_{i}$ of order $k_{i}$ is the starting point of exactly $k_{i}$ horizontal separatrices.
\item
An area form $dw_{X} \wedge d\bar{w_{X}}$ on $X\setminus\Sigma$.
\end{itemize}

We recall from the introduction the notion of \emph{saddle connection} for a translation surface $X$. It is a geodesic path $\ga:(0,T)\to X$ for the flat metric such that $\ga^{-1}(\Sigma)=\{0,T\}$, that is it starts end ends at points in $\Sigma$ and does not contain other such points in its interior.

\subsubsection{Moduli space}\label{modulispace(sss)}

Let us consider two elements $X$ and $X'$ in $\trasl(k_{1},..,k_{r})$ whose sets of singular points are respectively $\Sigma=\{p_{1},..,p_{r}\}$ and $\Sigma'=\{p'_{1},..,p'_{r}\}$. We define an equivalence relation on $\trasl(k_{1},..,k_{r})$ saying that $X\sim X'$ if there exist a diffeomorphism $f:X\to X'$ such that
\begin{itemize}
\item
$f(p_{i})=p'_{i}$ for all  $i=1,..,r$
\item
$\rho$ is a translation chart of $X$ if and only if $\rho\circ f$ is a translation chart of $X'$.
\end{itemize}
The quotient is the \emph{stratum} of the so called \emph{moduli space} of translation surfaces with cone singularities of orders $k_{1},..,k_{r}$, and is denoted $\cH(k_{1},..,k_{r})$.

\subsubsection{Teichm\"uller space}\label{teichmullerspace(sss)}

We fix a reference topological pair $(M^{\ast},\Sigma^{\ast})$, where $M^{\ast}$ is homeomorphic to any surface in $\trasl(k_{1},..,k_{r})$ and $\Sigma^{\ast}=\{p^{\ast}_{1},..,p^{\ast}_{r}\}$ is any subset of $M^{\ast}$ with $r$ elements. We denote $\diff^{+}(M^{\ast},\Sigma^{\ast})$ the group of diffeomorphisms of $M^{\ast}$ which preserve the orientation and are the identity on $\Sigma^{\ast}$ and we call $\diff^{+}_{0}(M^{\ast},\Sigma^{\ast})$ the subgroup of those diffeomorphisms which are isotopic to the identity. A \emph{marking} of a translation surface $X$ is a diffeomorphism of pairs
$$
\phi:(M^{\ast},\Sigma^{\ast})\to (X,\Sigma)
$$
such that $\phi(p^{\ast}_{i})=p_{i}$ for all $i=1,..,r$. We denote $(X,\phi)$ the datum of a translation surface with marking and we denote $\trasl^{\ast}(k_{1},..,k_{r})$ the set of translation surfaces with marking. We define an equivalence relation saying that $(X,\phi)\sim (X',\phi')$ if there exist a diffeomorphism $f:X\to X'$ as in paragraph \ref{modulispace(sss)} such that
$$
\phi'^{-1}\circ f\circ \phi \in \diff^{+}_{0}(M^{\ast},\Sigma^{\ast}).
$$
The quotient of $\trasl^{\ast}(k_{1},..,k_{r})$ by $\sim$ is a stratum of the so-called \emph{Teichm\"uller space} of translation surfaces with cone singularities of orders $k_{1},..,k_{r}$, and is denoted $\cT(k_{1},..,k_{r})$.

\subsubsection{Period map and local coordinates}\label{periodmap(sss)}

We fix a basis of cycles $\{\ga^{\ast}_{1},..,\ga^{\ast}_{d}\}$ for $H_{1}(M^{\ast},\Sigma^{\ast},\ZZ)$, where $d=2g+r-1$. The so called \emph{period map} $\Th:\cT(k_{1},..,k_{r})\to \CC^{d}$ associates to a translation structure with marking $(X,\phi)$ a complex vector $z$ in $\CC^{d}$ by the formula
\begin{equation}\label{eqperiodmap}
z_{i}:=\int_{\phi(\ga^{\ast}_{i})}w_{X}.
\end{equation}
We observe that the integral in (\ref{eqperiodmap}) is an isotopy invariant, so $\Th$ is well defined. The Teichm\"uller space $\cT(k_{1},..,k_{r})$ inherits its topology from $\CC^{d}$ via $\Th$.

\begin{propnonnum}
The map $\Th$ in (\ref{eqperiodmap}) is a local homeomorphism.
\end{propnonnum}

\begin{demo}
The map $\Th$ is continuous by definition, so it is enough to prove that it is open and that its restriction to small open sets is injective. It is known that any translation surface $X$ admits a triangulation whose vertices are the points in $\Sigma$ and whose edges are saddle connections (see \cite{esma} for example). If $(X,\phi)$ is a translation surface with marking, the pre-images under $\phi$ of the edges of the triangulation are elements in $H_{1}(M^{\ast},\Sigma^{\ast},\ZZ)$ and the triangulation can be represented combinatorially as a (non planar) graph, whose vertices are the points in $\Sigma^{\ast}$ and whose edges are elements in $H_{1}(M^{\ast},\Sigma^{\ast},\ZZ)$. Once the combinatorics of the triangulation is fixed, the values $z_{1},..,z_{d}$ of the periods on the basis $\{\ga^{\ast}_{1},..,\ga^{\ast}_{d}\}$ of $H_{1}(M^{\ast},\Sigma^{\ast},\ZZ)$ determine the triangulation, and therefore the translation surface. In particular the values $z_{1},..,z_{d}$ for the periods which define a triangulation are determined by some system of linear inequalities, which give an open condition in $\CC^{d}$, therefore the map $\Th$ is open. On the other hand, we can take open sets $U$ in $\cT(k_{1},..,k_{r})$ small enough in order to have that any $(X,\phi)$ in $U$ admits a triangulation whose combinatorial representation is the same. From the discussion above it follows that any such $(X,\phi)$ is determined by the values $z_{i}=\int_{\phi(\ga^{\ast}_{i})}w_{X}$ of its periods, therefore $\Th|_{U}$ is injective. The lemma is proved.
\end{demo}

Since $\Th$ is a local homeomorphism, $\cT(k_{1},..,k_{r})$ inherits from $\CC^{2g+r-1}$ the structure of complex manifold, whose complex dimension is of course $2g+r-1$. The group $\diff^{+}(M^{\ast},\Sigma^{\ast})$ acts (on the right) on $\cT(k_{1},..,k_{r})$ by
$$
(X,\phi)\mapsto (X,\phi\circ f),
$$
moreover the action of the subgroup $\diff^{+}_{0}(M^{\ast},\Sigma^{\ast})$ is trivial, since the Teichm\"uller space is defined modulo isotopies. It follows that we have an action (on the right) on $\cT(k_{1},..,k_{r})$ of the so called \emph{mapping class group}
$$
\modu(g,r):=\diff^{+}(M^{\ast},\Sigma^{\ast})/
\diff^{+}_{0}(M^{\ast},\Sigma^{\ast}).
$$
The action if faithful and proper, but it has non trivial stabilizer at surfaces admitting complex automorphisms. The quotient is exactly the moduli space $\cH(k_{1},..,k_{r})$, which inherits therefore the structure of a complex orbifold with complex dimension $2g+r-1$.

\subsubsection{Smooth measure and Teichm\"uller flow}\label{smoothmeasure&teichmullerflow(sss)}

For any $G\in\glduer$ and any $X\in\trasl(k_{1},..,k_{r})$ we define a new translation surface $GX$ whose local charts have the form $G\circ \varphi$, where $\varphi$ is a local chart of $X$. Since the translation subgroup is normal in the group of affine automorphisms of $\RR^{2}$, then $GX$ is still a translation surface. Thus $\glduer$ acts on the left on $\trasl(k_{1},..,k_{r})$. Since $\diff^{+}(M^{\ast},\Sigma^{\ast})$ acts on the right on $\trasl(k_{1},..,k_{r})$, then the two actions commute and we get an action (on the left) of $\glduer$ on $\cH(k_{1},..,k_{r})$

Let us consider the standard volume form $d\leb$ on $\CC^{2g+r-1}$, normalized in order to give co-volume one to the integer lattice $(\ZZ\oplus i\ZZ)^{2g+r-1}$. Using the period map $\Th$ we pull back the standard form and we get a smooth form $dm:=\Th^{\ast}d\leb$ on the stratum $\cT(k_{1},..,k_{r})$, whose associated smooth measure is denoted $m$. The action of any element $f$ in $\modu(g,r)$ on $H_{1}(M^{\ast},\Sigma^{\ast},\ZZ)$ is represented by some $A_{f}$ in $\sldz$ which satisfies $\Th\circ f=A_{f}\circ \Th$, thus $\modu(g,r)$ preserves the form $dm$ on $\cT(k_{1},..,k_{r})$, and therefore the measure $m$. It follows that the projection to the quotient induces a well defined volume element $d\mu$ on $\cH(k_{1},..,k_{r})$, whose associated smooth measure is denoted $\mu$.

For any $X\in \trasl(k_{1},..,k_{r})$ the corresponding holomorphic $1$-form $w_{X}$ induces an area form $w_{X}\wedge \overline{w}_{X}$ on $X\setminus \Sigma$. The area of $X$ is given by
$$
\area(X):=\int _{X\setminus \Sigma}w_{X}\wedge \overline{w}_{X}.
$$
The integral above is obviously invariant under the action of $\diff^{+}(M^{\ast},\Sigma^{\ast})$ and in particular of $\diff^{+}_{0}(M^{\ast},\Sigma^{\ast})$, therefore it defines a function on $\cT(k_{1},..,k_{r})$ which is invariant under $\modu(g,r)$. This amounts to say that we have an \emph{area function} $X\mapsto \area(X)$ on the moduli space $\cH(k_{1},..,k_{r})$. For any choice of a marking $(X,\phi)$ in $\cT(k_{1},..,k_{r})$ the well known \emph{Riemann's bilinear relation} expresses $\area(X)$ as a real analytic function of $z=\Th(X,\phi)$, it follows that the area function is smooth on $\cT(k_{1},..,k_{r})$ and we have a sub-manifold of real codimension one
$$
\cT^{(1)}(k_{1},..,k_{r}):=\{(X,\phi)\in\cT^{(1)}(k_{1},..,k_{r});\area(X)=1\}.
$$
The muduli space of area one translation surfaces is the smooth hyperboloid
$$
\cH^{(1)}(k_{1},..,k_{r}):=\cT^{(1)}(k_{1},..,k_{r})/\modu(g,r),
$$
it is possible to see that it is non-compact and in general non-connected (see \cite{konzo}). We have an homeomorphism from $\RR_{+}\times \cH^{(1)}(k_{1},..,k_{r})$ to $\cH(k_{1},..,k_{r})$ sending $(\la,X)$ to $\la X$, where the multiplication by $\la$ is given by the action of $\glduer$. It follows that the volume element $d\mu$ decomposes as
$$
d\mu=\la^{n}d\la\wedge d\mu^{(1)},
$$
where $n=4g+2r-3$ and $d\mu^{(1)}$ is a volume form on $\cH^{(1)}(k_{1},..,k_{r})$, whose associated smooth measure is denoted $\mu^{(1)}$.

$\slduer$ acts on $\cH^{(1)}(k_{1},..,k_{r})$ as subgroup of $\glduer$ and it is easy to see that it preserves $\mu^{(1)}$. The \emph{Teichm\"uller flow} $\cF_{t}$ on $\cH^{(1)}(k_{1},..,k_{r})$ is the action of the diagonal subgroup of $\slduer$, that if for any $t\in\RR$:
$$
\cF_{t}=
\left(
\begin{array}{cc}
e^{t} & 0 \\
0 & e^{-t}
\end{array}
\right).
$$
We recall the following fundamental basic result in Teichm\"uller dynamics (see \cite{masuruno} and \cite{veech}).
\begin{thmnonnum}
The smooth measure $\mu^{(1)}$ gives to $\cH^{(1)}(k_{1},..,k_{r})$ finite volume and its restriction on any connected component is ergodic with respect to $\cF_{t}$.
\end{thmnonnum}

\subsection{Interval exchange transformations and Veech's construction.}\label{ietandveechconstruction(ss)}

\subsubsection{Interval exchange transformations}\label{iet(sss)}

Let $X$ be a translation surface and let $I$ be an horizontal segment in $X$, that is a finite segment of a trajectory of the horizontal field $\partial_{x}$. For $x\in I$ we denote $T(x)$ the first intersection of $I$ with the positive trajectory of the vertical field $\partial_{y}$ starting from $x$, in other words we consider the \emph{first return map} to $I$ of the vertical flow. Since $\partial_{y}$ preserves the area element on $X$, then $T$ is defined almost everywhere on $I$, moreover $T$ preserves the length element $dx$ on $I$, therefore it acts locally as a translation in the coordinate $x$. It follows that the domain of $T$ is union of open sub-intervals of $I$, whose endpoints lie on positive trajectories of $\partial_{y}$ ending in points of $\Sigma$ without passing from $I$. There are just finitely many such trajectories, thus the domain of $T$ is union of a finite number of intervals and the effect of $T$ is to rearrange them in $I$ by translations. A map like $T$ is called an \emph{interval exchange transformation}.

\begin{defi}
Let $d\geq 2$ and let $\la=(\la_{1},..,\la_{d})$ be a vector in $\RR^{d}_{+}$ and $\pi$ be a permutation in the symmetric group $S_{d}$. An interval exchange transformation (we will write i.e.t.) with $d$ intervals is a map $T$ from an interval $I$ to itself defined by the data $(\pi,\la)$ as follows.
\begin{itemize}
\item
The interval $I$ admits two partitions into open sub-intervals $I=\sqcup_{i=1}^{d}I^{t}_{i}$ and $I=\sqcup_{i=1}^{d}I^{b}_{i}$. For any $i\in\{1,..,d\}$, if we start counting from the left, $I^{t}_{i}$ is in $i$-th position in the first partition and $I^{b}_{i}$ is in $i$-th position in the second partition.
\item
For any $i\in\{1,..,d\}$ we have that $|I^{t}_{i}|=|I^{b}_{\pi(i)}|=\la_{i}$ and the restriction of $T$ to $I^{t}_{i}$ is the translation onto $I^{b}_{\pi(i)}$.
\end{itemize}
The data $\pi$ and $\la$ are called respectively the \emph{combinatorial datum} and the \emph{length datum} of $T$.
\end{defi}

For any combinatorial datum $\pi$ let us call $\De_{\pi}:=\{\pi\}\times \RR^{d}_{+}$ the set of all i.e.t.s with combinatorial datum $\pi$. For an i.e.t. $T$ with $d$ intervals and for any $i\in\{1,..,d\}$ we call $u^{t}_{i}$ the left endpoint of $I^{t}_{i}$ and $u^{b}_{i}$ the left endpoint of $I^{b}_{i}$. Any $u^{t}_{i}$ is a discontinuity of $T$ and any $u^{b}_{i}$ is a discontinuity of $T^{-1}$. We say that the combinatorial datum $\pi$ is \emph{admissible} if for any integer $k$ with $1\leq k<d$ we have $\pi\{1,..,k\}\not =\{1,..,k\}$. A \emph{connection} for $T:I\to I$ is a triple $(q,p,n)$ with $1<q\leq d$, $1<p\leq d$ and $n\in \NN$ such that $T^{n}u^{b}_{q}=u^{t}_{p}$. In particular, if $T$ has no connection, then its combinatorial datum is admissible.

\subsubsection{The Veech construction.}\label{veechconstruction(sss)}

In this paragraph we describe a construction given by Veech in \cite{veech} and known as \emph{zippered rectangles construction}. We follow the presentation of \cite{mmy}, even if we use a slightly different notation. Let $d\geq 2$. For an admissible combinatorial datum $\pi$ in $S_{d}$ let $\Theta_{\pi}$ be the open convex polyhedral cone in $\RR^{d}$ of those $\tau$ such that for any $1\leq k \leq d-1$ we have
$$
\sum_{i\leq k} \tau_{i}>0 \texttt{ and } \sum_{\pi(i)\leq k} \tau_{\pi(i)}<0.
$$
The cone $\Theta_{\pi}$ is never empty since it contains the vector $\tau$ with coordinates $\tau_{i}:=\pi(i)-i$. We say that a vector $\tau$ in $\Theta_{\pi}$ is a \emph{suspension datum} for $\pi$.

Let us consider a pair $(\la,\tau)$ with $\la\in\De_{\pi}$ and $\tau\in\Th_{\pi}$. The data $(\pi,\la)$ define an i.e.t. $T$ on some interval $I$. We will construct a translation surface $X=X(\pi, \la, \tau)$ containing $I$ as horizontal segment and such that the construction in paragraph \ref{iet(sss)} gives exactly $T$ as first return to $I$ for the vertical field of $X$. We first define the complex vector $\zeta =\lambda +\sqrt{-1}\tau \in \CC^{d}$, then for any $i\in\{1,..,d\}$ we define two complex numbers
$$
\xi_{i}^{t}:=\sum_{j<i}\zeta_{j}
\texttt{ and }
\xi_{i}^{b}:=\sum_{\pi(j)<\pi(i)}\zeta_{j}.
$$
We observe that $\xi^{t}_{1}=\xi^{b}_{1}=0$ and we set $\xi_{\ast}=\sum_{i=1}^{d}\zeta_{i}$, $\la_{\ast}:=\sum_{i=1}^{d}\la_{i}$ and $\tau_{\ast}=\sum_{i=1}^{d}\tau_{i}$. Condition $\tau\in\Th_{\pi}$ means that for any $i$ with $1<i\leq d$ we have $\Im(\xi_{i}^{t})>0$ and $\Im(\xi_{i}^{b})<0$. We also define the \emph{translation datum} $\theta \in \CC^{d}$ setting $\theta_{i}:=\xi^{b}_{\pi(i)}-\xi^{t}_{i}$ for any $i$, then we decompose it as $\theta = \de-\sqrt{-1}h$ with $\de,h\in\RR^{d}$, in particular, since $\pi$ is admissible, we have $h_{i}=\Im (\xi^{t}_{i})-\Im (\xi^{b}_{\pi(i)})>0$ for any $i\in\{1,..,d\}$.

We may suppose that the i.e.t. $T$ defined by the data $(\pi,\la)$ acts on the interval $I=(0,\la_{\ast})$. It $u^{t}_{i}$ and $u^{b}_{i}$ denote the singularities for $T$ and $T^{-1}$ we have
$$
u^{t}_{i}:=\sum_{j<i}\lambda_{j}
\texttt{ and }
u^{b}_{i}:=\sum_{\pi(j)<\pi(i)}\lambda_{j}.
$$
Observe that $u_{i}^{t}=\Re(\xi_{i}^{t})$ and $u^{b}_{i}=\Re(\xi_{i}^{b})$ for all $i=1,..,d$. We embed $I$ in the complex plane, that is we consider $I=(0,\la_{\ast})\times \{0\}$, then we define $2d$ open rectangles in $\CC$ setting for any $i\in\{1,..,d\}$
$$
\begin{array}{c}
R_{i}^{t}:=(u^{t}_{i},u^{t}_{i}+\lambda_{i})\times (0,h_{i})\\
R_{i}^{b}:=(u^{b}_{i},u^{b}_{i}+\lambda_{\pi^{-1}(i)})\times (-h_{\pi^{-1}(i)},0).
\end{array}
$$
In order to get a surface we zip together these rectangles with the identifications described below.

\begin{enumerate}
\item
For each $i\in\{1,..,d\}$ the rectangle $R_{i}^{t}$ is equivalent to the rectangle $R_{\pi(i)}^{b}$ under the translation by the translation datum $\te_{i}$.
\item
For each $i$ with $1<i\leq d$ we paste together $R^{t}_{i}$ and $R^{t}_{i-1}$ along the common vertical open segment which connects the point $u^{t}_{i}+0\sqrt{-1}$ to the point $\xi^{t}_{i}$.
\item
For each $i$ with $1<i\leq d$ we paste together $R^{b}_{i}$ and $R^{b}_{i-1}$ along the vertical open segment which connects the point $u^{b}_{i}+0\sqrt{-1}$ to the point $\xi^{b}_{i}$.
\item
For any $i\in\{1,..,d\}$ we paste $R_{i}^{t}$ to $I$ along its lower horizontal boundary segment $(u^{t}_{i},u^{t}_{i}+\la_{i})\times \{0\}$.
\item
For any $i\in\{1,..,d\}$ we paste $R_{i}^{b}$ to $I$ along its upper horizontal boundary segment $(u^{b}_{i},u^{b}_{i}+\la_{\pi^{-1}(i)})\times \{0\}$.
\item
We observe that $h_{\pi^{-1}(d)}= \Im(\xi^{t}_{\pi^{-1}(d)})-\Im(\xi^{b}_{d})= \Im(\xi^{t}_{\pi^{-1}(d)+1})-\tau_{\ast}$ and similarly $h_{d}=\tau_{\ast}-\Im(\xi^{b}_{\pi(d)+1})$. If $\tau_{\ast}\geq 0$ then $h_{\pi^{-1}(d)}\geq\Im(\xi^{t}_{\pi^{-1}(d)+1})$ and the translation by $\te_{\pi^{-1}(d)}$ induces an identification between the vertical segment connecting $u^{t}_{\pi^{-1}(d)+1}+\sqrt{-1}h_{\pi^{-1}(d)}$ to $\xi^{t}_{\pi^{-1}(d)+1}$ and the vertical segment connecting $\la_{\ast}+0\sqrt{-1}$ to $\xi_{\ast}$. If $\tau_{\ast}<0$ then $-h_{d}>\Im(\xi^{b}_{\pi(d)+1})$ and the translation by $-\te_{d}$ induces an identification between the vertical segment connecting $\xi^{b}_{\pi(d)+1}$ to $u^{b}_{\pi(d)+1}-\sqrt{-1}h_{d}$ and the vertical segment connecting $\xi_{\ast}$ to $\la_{\ast}+0\sqrt{-1}$.
\item
Finally we add the origin $0$ of $\CC$, the points $\xi^{t}_{i}$ with $1<i\leq d$, the points $\xi^{b}_{i}$ with $1<\pi(i)\leq d$ and the point $\xi_{\ast}$
%(note that $\xi_{\al'_{t}}^{t}=\xi_{\al'_{b}}^{b}$ if $\al'_{t}$ and $\al'_{b}$ are such that $\pi^{t}(\al'_{t})=\pi^{b}(\al'_{b})=1$)
\end{enumerate}

\begin{defi}\label{def1theveechconstruction}
For any admissible combinatorial datum $\pi$ in $S_{d}$ and for any pair of length-suspension data $(\la,\tau)$ for $\pi$ we call $X(\pi,\la,\tau)$ the translation surface which is obtained following the procedure above. A translation surface $X$ such that there exists a triple of data $(\pi,\la,\tau)$ as above satisfying $X=X(\pi,\la,\tau)$ is said \emph{representable with the Veech construction}.
\end{defi}

Now we determine the stratum of the moduli space where the surface $X(\pi,\la,\tau)$ lies. We fix a small positive number $\epsilon$ and for $i\in\{2,..,d\}$ we define two half sectors $D^{t}_{i}:=B(\xi^{t}_{i},\epsilon)\cap R^{t}_{i}$ and $D^{b}_{i}:=B(\xi^{b}_{i},\epsilon)\cap R^{b}_{i}$ in the complex plane,
where the notation suggests that these sectors look like the letter D. Any $D^{t}_{i}$ is equivalent to $D^{b}_{\pi(i)}$ under the translation by $\te_{i}$, therefore they are identified to the same sector in $X(\pi,\la,\tau)$. Similarly for $i\in\{1,..,d-1\}$ we define two half sectors $C^{t}_{i}:=B(\xi^{t}_{i},\epsilon)\cap R^{t}_{i+1}$ and $C^{b}_{i}:=B(\xi^{b}_{i},\epsilon)\cap R^{b}_{i+1}$, here the notation suggests that these sectors look like the letter C. As before $C^{t}_{i}$ is equivalent to $C^{b}_{\pi(i)}$ under the translation by $\te_{i}$ and they are identified to the same sector in $X(\pi,\la,\tau)$. We introduce the set $HS:=\{(i,D)\}_{i\in\{2,..,d\}}\cup\{(i,C)\}_{i\in\{1,..,d-1\}}$, which has $2(d-1)$ elements, that we use as labels for the sectors introduced above. The labeling obviously associates the symbol $(i,D)$ to any half sector $D^{t}_{i}$ and the symbol $(i,C)$ to any half sector $C^{t}_{i}$. Finally we introduce a bijection $\widehat{\pi}$ of the elements of $HS$ defined by
$$
\widehat{\pi}(i,D):=(\pi(i-1),C)
\texttt{ and }
\widehat{\pi}(i,C):=(\pi^{-1}(i+1),D).
$$
Let us consider any singular point $p\in\Sigma$ for the surface $X(\pi,\la,\tau)$ and any half sector $D^{t}_{i}$ in $p$. Turning in counterclockwise sense around $p$ we pass from $D^{t}_{i}$ to $C^{t}_{i-1}$, which is identified to $C^{b}_{\pi(i-1)}$. Then turning in the same sense we pass from $C^{b}_{\pi(i-1)}$ to $D^{b}_{\pi(i-1)+1}$, which is equivalent to $D^{t}_{\pi^{-1}(j+1)}$, where $j=\pi(i-1)$. It follows that for any $p\in \Sigma$, the number ho half sectors that we meet turning around $p$ equals to the length of a cycle of $\widehat{\pi}$. We have proved the following

\begin{lemnonnum}
For any triple of data $(\pi,\la,\tau)$ the translation surface $X(\pi,\la,\tau)$ lies in a stratum $\cH(k_{1},..,k_{r})$ which depends only from the combinatorial datum $\pi$.
\end{lemnonnum}

We also recall the following well known result (see \cite{ytre}).

\begin{lemnonnum}
For any $X$ in $\cH(k_{1},..,k_{r})$ without vertical saddle connections there exists an admissible $\pi\in S_{2g+r-1}$ and a pair of data $(\la,\tau)\in\De_{\pi}\times\Th_{\pi}$ such that $X=X(\pi,\la,\tau)$, anyway the triple of data $(\pi,\la,\tau)$ representing $X$ is not unique. In particular almost any $X$ in $\cH(k_{1},..,k_{r})$ is representable with the Veech construction.
\end{lemnonnum}

\subsubsection{Veech local charts.}\label{veechlocalcharts(sss)}

We recall that in the definition of the Teichm\"uller space we considered a reference topological pair $(M^{\ast},\Sigma^{\ast})$, moreover in order to define the period map we fixed also a basis $\{\ga^{\ast}_{1},..,\ga^{\ast}_{d}\}$ for $H_{1}(M^{\ast},\Sigma^{\ast},\ZZ)$.

According to the construction in paragraph \ref{veechconstruction(sss)}, for $d\geq 2$ and for any admissible $\pi\in S_{d}$ we have a well defined map
\begin{equation}\label{eq1veechlocalcharts}
\begin{array}{cccc}
\cI_{\pi}: & \De_{\pi}\times \Th_{\pi} &   \to   & \cH(k_{1},..,k_{r})\\
           &    (\la,\tau)             & \mapsto & X(\pi,\la,\tau).
\end{array}
\end{equation}
Let us consider any $(\la,\tau)\in \De_{\pi}\times \Th_{\pi}$ and the associated translation surface $X=X(\pi,\la,\tau)$. Using the notation of paragraph \ref{veechconstruction(sss)}, for any $i\in\{1,..,d\}$ we define $\widehat{\zeta}^{t}_{i}$ as the segment in $\CC$ connecting $\xi^{t}_{i}$ to $\xi^{t}_{i}+\zeta_{i}$ and $\widehat{\zeta}^{b}_{i}$ as the segment in $\CC$ connecting $\xi^{b}_{i}$ to $\xi^{b}_{i}+\zeta_{\pi(i)}$. It is easy to see that for any $i\in\{1,..,d\}$, either $\widehat{\zeta}^{t}_{i}$ or $\widehat{\zeta}^{b}_{i}$ projects to a saddle connection in $X(\pi,\la,\tau)$, that we denote $\widehat{\zeta}_{i}$. It is also easy to see that the set of curves $\{\widehat{\zeta}_{1},..,\widehat{\zeta}_{d}\}$ is a basis for $H_{1}(X,\Sigma,\ZZ)$. Therefore we can define naturally an isotopy class $\phi(\pi,\la,\tau)$ of maps of pairs $\phi:(M^{\ast},\Sigma^{\ast})\to (X,\Sigma)$ such that for any $i\in\{1,..,d\}$ we have
$$
\phi(\pi,\la,\tau)(\ga^{\ast}_{i})=\widehat{\zeta}_{i},
$$
where the equal means that the two curves are \emph{isotopic} (non just \emph{homologous}). It follows that the zippered rectangle construction induces naturally a marking $\phi(\pi,\la,\tau)$ of $X(\pi,\la,\tau)$, that is we have a map

\begin{equation}\label{eq2veechlocalcharts}
\begin{array}{cccc}
\cI^{\cT}_{\pi}: & \De_{\pi}\times \Th_{\pi} &   \to   & \cT(k_{1},..,k_{r})\\
           &    (\la,\tau)             & \mapsto & (X(\pi,\la,\tau),\phi(\pi,\la,\tau)).
\end{array}
\end{equation}

\begin{lem}\label{lemveechlocalcharts}
For any admissible $\pi$ the map $\cI^{\cT}_{\pi}$ in (\ref{eq2veechlocalcharts}) is an homeomorphism onto its image, moreover the push-forward $\cI^{\cT}_{\pi\ast}\leb$ of the lebesge measure on $\De_{\pi}\times \Th_{\pi}$ is equivalent to the restriction of the smooth measure $m$ to the image of $\cI^{\cT}_{\pi}$.
\end{lem}

\begin{demo}
Since $\phi(\pi,\la,\tau)(\ga_{i}^{\ast})=\widehat{\zeta}_{i}$ for any $i\in\{1,..,d\}$, then it follows from equation (\ref{eqperiodmap}) that $\Th \circ \cI^{\cT}_{\pi}$ is the natural immersion of $\De_{\pi}\times \Th_{\pi}$ in $\CC^{2g+r-1}$, where $g$ and $r$ are determined by $\pi$. Since $\Th$ is a local homeomorphism, then also $\cI^{\cT}_{\pi}$ is, moreover it is injective, therefore is an homeomorphism. Since $\cI^{\cT}_{\pi}$ is a local inverse of $\Th$ the second part of the statement is obvious.
\end{demo}

\begin{cor}\label{corveechlocalcharts}
For any admissible $\pi$ the map $\cI_{\pi}$ in (\ref{eq1veechlocalcharts}) is continuous and open, moreover the push-forward $\cI_{\pi\ast}\leb$ of the lebesge measure on $\De_{\pi}\times \Th_{\pi}$ is equivalent to the restriction of the smooth measure $\mu$ to the image of $\cI_{\pi}$.
\end{cor}

\section{Dichotomy for the Khinchin type condition}

%This section is devoted to the proof of theorem \ref{teoremab}.

\subsubsection{A remark on strongly decreasing functions}\label{remarkstronglydecreasingfunctions}

We recall from the introduction the notion of strongly decreasing function, that is $\varphi:[0,+\infty)\to(0,\infty)$ such that $t\varphi(t)$ is decreasing monotone. We can restrict such function to integers, obtaining a sequence $(\varphi(n))_{n\in\NN}$. In particular $\varphi(t)$ is monotone, therefore $\sum_{n=0}^{\infty}\varphi(n)$ is of the same order of $\int_{0}^{+\infty}\varphi(t)dt$, that is the series diverges if and only if the integral diverges. For any real number $\te>1$ we can also consider the sequence $\psi_{n}:=\te^{n}\varphi(\te^{n})$ and it is an easy exercise in calculus to see that since $t\varphi(t)$ is decreasing monotone than $\sum_{n=0}^{\infty}\psi_{n}=+\infty$ if and only if $\int_{0}^{+\infty}\varphi(t)dt=+\infty$.

\subsection{Proof of the convergent case of theorem \ref{teoremab}}\label{prooftheorembconvergentcase(ss)}

In this paragraph we prove proposition \ref{proposizioneb} and therefore the convergent part of theorem \ref{teoremab} too. Let $X$ be any translation surface. $\soduer$ acts on the stratum of $X$ as subgroup of $\slduer$ and we call $X_{\te}$ the image of $X$ under
$$
R_{\te}=
\left(
\begin{array}{cc}
  \cos\te & -\sin\te \\
  \sin\te & \cos\te
\end{array}
\right).
$$
Let $\varphi:[0,+\infty)\to (0,+\infty)$ be a strongly decreasing function with $\int_{0}^{+\infty}\varphi(t)dt<+\infty$. We show that for almost any $\te\in[0,2\pi)$ the set $\hol(X_{\te})$ contains just finitely many solutions of equation (\ref{eqphiconnessionesuptrasl}). Since the flat metric of $X_{\te}$ is the same as the flat metric of $X$, then the set of all the saddle connections of $X$ is also the set of the saddle connections of $X_{\te}$ for any $\te$. It follows that $\hol(X_{\te})=R_{\te}\hol(X)$ and in particular the length $|\hol(\ga)|$ of a saddle connection $\ga$ is invariant under rotations. For any $v\in\hol(X)$ we define the set
$$
I(v):=\{\te\in [0,2\pi);|\Re(R_{\te}v)|\leq \varphi (|v|)\},
$$
that is the set of those $\te$ such that $\hol(X_{\te})$ contains $R_{\te}v$ as solution of equation (\ref{eqphiconnessionesuptrasl}). Our strategy is to prove that $\sum_{v\in\hol(X)}\leb (I(v))<\infty$, then the classical Borel-Cantelli lemma implies that almost any $\te\in [0,2\pi)$ is contained into a finite number of sets $I(v)$, which is equivalent to say that $\hol(X_{\te})$ contains just finitely many solutions of equation (\ref{eqphiconnessionesuptrasl}) for almost any $\te$.

In \cite{masurdue} Masur proves that for any translation surface $X$ the number $N(X,L)$ of elements $v$ of $\hol(X)$ with norm $|v|\leq L$ has quadratic growth with $L$, that is there are two positive constants $c<C$ such that for any $L$ big enough we have
$$
cL^{2}\leq N(X,L)\leq CL^{2}.
$$
We write
$$
\sum_{v\in\hol(X)}\leb (I(v))=
\sum_{|v|\leq 1}\leb (I(v))+
\sum_{k\in \NN^{\ast}}
\left(\sum_{2^{k-1}<|v|\leq 2^{k}}\leb (I(v))\right).
$$
We observe that $\Re(R_{\te}v)/|v|$ is the sinus of the angle between $R_{\te}v$ and the imaginary axis. Since $\varphi$ is bounded, than for any $\epsilon >0$ and for $|v|$ big enough we have $\leb (I(v))<(2+\epsilon)\varphi(|v|)/|v|$. Recalling that $t\varphi(t)$ is decreasing monotone, and therefore $\varphi(t)$ too, and applying Masur's estimation for $N(X,2^{k})$ with $k$ big enough, we get
$$
\sum_{2^{k-1}<|v|\leq 2^{k}}\leb (I_{\ga})\leq C(2+\epsilon)2^{2k}\frac{\varphi(2^{k-1})}{2^{k-1}}.
$$
It follows that the tail of the sum $\sum_{v\in\hol(X)}\leb (I(v))$ is controlled by the series $\sum_{k=1}^{\infty}2^{k}\varphi(2^{k})$, modulo a multiplicative constant. According to remark \ref{remarkstronglydecreasingfunctions}, condition $\int_{0}^{+\infty}\varphi(t)dt<+\infty$ implies $\sum_{v\in\hol(X)}\leb (I(v))<+\infty$. Proposition \ref{proposizioneb} is proved and the convergent part of theorem \ref{teoremab} too.

\subsection{Bundles of saddle connections}\label{bundlesofsaddleconnections(ss)}

\subsubsection{Framed translation surfaces}\label{framedtranslationsurfaces(sss)}

A \emph{frame} for a translation surface $X$ is the datum of $r$ horizontal separatrices $(S_{1},..,S_{r})$ such that $S_{i}$ starts in $p_{i}\in\Sigma$ for any $i=1,..,r$. We denote $\widehat{X}$ the datum $(X,S_{1},..,S_{r})$ of a translation surface with frame and $\widehat{\trasl}(k_{1},..,k_{r})$ the set of translation surfaces with frame whose singularities have orders $k_{1},..,k_{r}$. An equivalence relation is defined on this set setting $\widehat{X}\sim \widehat{X}'$ if there exists an diffeomorphism $f:X\to X'$ as in the definition of the moduli space in paragraph \ref{modulispace(sss)} such that for any $i=1,..,r$ we have
$$
f(S_{i})=S'_{i},
$$
where $(S_{1},..,S_{r})$ and $(S'_{1},..,S'_{r})$ are respectively the frame of $\widehat{X}$ and the frame of $\widehat{X}'$. We denote $\widehat{\cH}(k_{1},..,k_{r})$ the quotient, which is called the moduli space of translation surfaces with frame. The map
$$
\begin{array}{ccc}
\widehat{\cH}(k_{1},..,k_{r}) &   \to   & \cH(k_{1},..,k_{r})\\
\widehat{X}  & \mapsto & X
\end{array}
$$
is a covering map of degree $\prod_{i=1}^{r}k_{i}$, therefore $\widehat{\cH}(k_{1},..,k_{r})$ inherits the structure of complex orbifold of complex dimension $2g+r-1$. Obviously all the structures defined on $\cH(k_{1},..,k_{r})$ have a corresponding lift in the covering space. In particular we have an action of $\glduer$ and a smooth measure on $\widehat{\cH}(k_{1},..,k_{r})$, that we still call $\mu$, which is invariant under the sub-group action of $\slduer$. The smooth hyperboloid $\widehat{\cH}^{(1)}(k_{1},..,k_{r})$ of area one translation surfaces with frame is also obviously defined, together with the induced smooth $\slduer$-invariant volume $\mu^{(1)}$. We will be interested to the action of the Teichm\"uller flow, that is the action of the diagonal subgroup
$$
\cF_{t}=
\left(
\begin{array}{cc}
e^{t} & 0 \\
0 & e^{-t} \\
\end{array}
\right).
$$
The definition of $\widehat{\cH}(k_{1},..,k_{r})$ is conceptually the same as the definition of the double-orientable cover of a non-orientable manifold. The introduction of the covering is necessary because the choice of a frame is not invariant under $\cF_{t}$ for elements of $\cH(k_{1},..,k_{r})$. For example in $\cH(k_{1},..,k_{r})$ there are closed orbits of $\cF_{t}$ with period $T$ such that, for any $X$ in the orbit and any choice of a frame $(S_{1},..,S_{r})$ on $X$, we have $(\cF_{T}S_{1},..,\cF_{T}S_{r})\not =(S_{1},..,S_{r})$.

\subsubsection{Bundles of saddle connections}\label{bundlesofsaddleconnections(sss)}

Let $\widehat{X}$ any element in $\widehat{\cH}(k_{1},..,k_{r})$, let $p_{i},p_{j}$ be any pair of point in the singular set $\Sigma$ of $X$ and let $l\in\{1,..,k_{i}\}$ and $m\in\{1,..,k_{j}\}$.

\begin{defi}\label{defconfigurazioneconnessionisella}
A \emph{bundle} of saddle connections is a set $\cC^{(p_{j},p_{i},m,l)}(\widehat{X})$ of those saddle connection $\ga$ for the translation surface $X$ which start in $p_{j}$, end in $p_{i}$ and satisfy
$$
\begin{array}{c}
2(m-1)\pi\leq \angolo(S_{j},\ga)<2m\pi\\
2(l-1)\pi\leq \angolo(S_{i},\ga)<2l\pi,
\end{array}
$$
where $(S_{1},..,S_{r})$ is the frame for $X$ corresponding to $\widehat{X}$. We denote $\hol^{(p_{j},p_{i},m,l)}(\widehat{X})$ the set of complex vectors $v=\int_{\ga}w_{X}$ for $\ga$ in $\cC^{(p_{j},p_{i},m,l)}(\widehat{X})$.
\end{defi}

\begin{lem}\label{lemconfigurazioneconnessionisella}
Let $\widehat{X}$ be any element in $\widehat{\cH}(k_{1},..,k_{r})$. Let $p_{j},p_{i}$ be any pair of points in $\Sigma$ and $(m,l)$ be a pair of integers with $1\leq m\leq k_{j}$ and $1\leq l\leq k_{i}$. Then, for any $\ga\in\cC^{(p_{j},p_{i},m,l)}(\widehat{X})$ and for any $t\in\RR$, $\ga$ is a saddle connection for $\cF_{t}X$ and moreover it belongs to $\cC^{(p_{j},p_{i},m,l)}(\cF_{t}\widehat{X})$.
\end{lem}

\begin{demo}
The Teichm\"uller flow obviously preserves saddle connections and the names of singular points in $\Sigma$. Let $p_{j}$ and $p_{i}$ be points in $\Sigma$ and let $\ga$ a saddle connection for $X$ starting at $p_{j}$ and ending at $p_{i}$. Let $(S_{1},..,S_{r})$ be the frame on $X$ carried by $\widehat{X}$. For any $t\in\RR$ we call $\ga_{t}$ the saddle connection for $\cF_{t}X$ corresponding to $\ga$. The frame carried by $\cF_{t}\widehat{X}$ is $(\cF_{t}S_{1},..,\cF_{t}S_{r})$ and we have
$$
\tan\angolo(\cF_{t}S_{j},\ga_{t})=
e^{-2t}\tan\angolo(S_{j},\ga),
$$
$$
\tan\angolo(\cF_{t}S_{i},\ga_{t})=
e^{-2t}\tan\angolo(S_{i},\ga),
$$
therefore the relations in definition \ref{defconfigurazioneconnessionisella} are preserved. The lemma is proved.
\end{demo}

\subsubsection{Relation with the Veech construction}\label{relationveech(sss)}

Let us fix $d\geq 2$ and an admissible combinatorial datum $\pi\in S_{d}$. For $(\la,\tau)$ in $\De_{\pi}\times\Th_{\pi}$ let us consider the translation surface $X(\pi,\la,\tau)$. We fix a small positive real number $\epsilon$ and for any $i\in\{2,..,d\}$ we define an horizontal segment $S^{\epsilon}_{i}:(0,\epsilon)\to \CC$ in the complex plane by
$$
S^{\epsilon}_{i}(t):=\xi^{t}_{i}+t(1,0).
$$
The projection to $X(\pi,\la,\tau)$ of the segment $S^{\epsilon}_{i}$ is the beginning segment of an horizontal separatrix in $X(\pi,\la,\tau)$, that we call $S^{Veech}_{i}$. There are exactly $d-1$ horizontal separatrices as above and they are all distinct. The orders of the singularities of $X(\pi,\la,\tau)$ satisfy $k_{1}+..+k_{r}=2g+r-2=d-1$, therefore
$$
\{S^{Veech}_{2},..,S^{Veech}_{d}\}
$$
is the set of all the horizontal separatrix for $X(\pi,\la,\tau)$.

\begin{defi}\label{defconfigurazioneconnessionisellaveech}
For a pair of integers $(q,p)$ with $1<q\leq d$ and $1<p\leq d$ we denote $\cV^{(q,p)}(\pi,\la,\tau)$ the set of those saddle connections $\ga$ for $X(\pi,\la,\tau)$ which start in the point of $\Sigma$ where $S^{Veech}_{q}$ starts, end in the point of $\Sigma$ where $S^{Veech}_{p}$ starts, and satisfy
$$
\begin{array}{c}
0\leq \angolo(S^{Veech}_{q},\ga)<2\pi\\
0\leq \angolo(S^{Veech}_{p},\ga)<2\pi.
\end{array}
$$
\end{defi}

\begin{lem}\label{lemconfigurazioneconnessionisellaveech}
Let $(\la,\tau)$ in $\De_{\pi}\times \Th_{\pi}$, consider $X=X(\pi,\la,\tau)$ in $\cH(k_{1},..,k_{r})$ and take any pre-image $\widehat{X}$ of $X$ in $\widehat{\cH}(k_{1},..,k_{r})$. Then for any pair of points $p_{j},p_{i}$ in $\Sigma$ and any pair of integers $(m,l)$ with $1\leq m\leq k_{j}$ and $1\leq l\leq k_{i}$ there exist an unique pair $(q,p)$ in $\{2,..,d\}^{2}$ such that
$$
\cC^{(p_{j},p_{i},m,l)}(\widehat{X})=\cV^{(q,p)}(\pi,\la,\tau).
$$
\end{lem}

\begin{demo}
Let $(S_{1},..,S_{r})$ be the frame on $X=X(\pi,\la,\tau)$ carried by $\widehat{X}$. Since $(S^{Veech}_{2},...,S^{Veech}_{d})$ are all the horizontal separatrices of $X$, then for any pair of points $p_{j},p_{i}$ in $\Sigma$ and any pair of integers $(m,l)$ with $1\leq m\leq k_{j}$ and $1\leq l\leq k_{i}$ there exist an unique pair $(q,p)$ in $\{2,..,d\}^{2}$ such that $S^{Veech}_{q}$ starts in $p_{j}$ (where $S_{j}$ starts), $S^{Veech}_{p}$ starts in $p_{i}$ (where $S_{i}$ starts) and we have
$$
\angolo(S^{Veech}_{q},S_{j})=2(m-1)\pi
\texttt{ and }
\angolo(S^{Veech}_{p},S_{i})=2(l-1)\pi.
$$
The statement in the lemma is therefore evident.
\end{demo}

\subsection{Combinatorially defined saddle connections}\label{combdefinedsaddleconnections(ss)}

Let $\pi\in S_{d}$ be an admissible combinatorial datum and $(\la,\tau)$ be a pair of length-suspension data in $\De_{\pi}\times \Th_{\pi}$. We recall that if $T:I\to I$ is the i.e.t. defined by the data $(\pi,\la)$, then the interval $I$ embeds naturally in the translation surface $X(\pi,\la,\tau)$ as an horizontal segment and $T$ is the first return map to $I$ of the vertical flow on $X(\pi,\la,\tau)$.

\subsubsection{Combinatorially defined homology classes}\label{combdefinedhomologyclasses(sss)}

In the following we use the notation of the Veech construction for $X=X(\pi,\la,\tau)$ introduced in paragraph \ref{veechconstruction(sss)}. We recall that to any complex vector $\zeta_{i}=\la+\sqrt{-1}\tau_{i}$ is naturally associated a saddle connection $\widehat{\zeta}_{i}$, which in particular defines an element $[\widehat{\zeta}_{i}]$ in $H_{1}(X,\Sigma,\ZZ)$. The formulae
$$
\xi_{i}^{t}=\sum_{j<i}\zeta_{j}
\texttt{   ;  }
\xi_{i}^{b}=\sum_{\pi(j)<\pi(i)}\zeta_{j}
\texttt{   ;  }
\te_{i}=\xi_{\pi(i)}^{b}-\xi_{i}^{t}
$$
extend formally on relative (or absolute) homology classes
$$
[\xi_{i}^{t}]:=\sum_{j<i}[\widehat{\zeta}_{i}]
\texttt{  ;   }
[\xi_{i}^{b}]:=\sum_{\pi(j)<\pi(i)}[\widehat{\zeta}_{i}]
\texttt{  ;   }
[\te_{i}]:=[\xi_{\pi(i)}^{b}]-[\xi_{i}^{t}].
$$
The homology classes above have a representant which is concatenation of saddle connections for the flat structure $X=X(\pi,\la,\tau)$, anyway in general it is not possible to find a representant which is a saddle connection itself.

We define a piecewise constant map $[\te]:I\to H_{1}(X,\ZZ)$ setting $[\te](x)=[\te_{i}]$ if $x\in I^{t}_{i}$, then we consider the Birkhoff sum $S_{n}[\te]$ of the function $[\te]$ over $T$, that is
$$
S_{n}[\te](x)=[\te](x)+..+[\te](T^{n-1}x).
$$
We fix a pair of indexes $q,p$ with $1<q\leq d$ and $1<p\leq d$. If $T$ has no connection we can iterate $T$ on $u^{b}_{q}$ infinitely many times and we get a sequence of elements in  $H_{1}(X,\Sigma,\ZZ)$ defined by
$$
[\ga]_{q,p,n}:=[\xi_{p}^{t}]-[\xi_{q}^{b}]-S_{n}[\te](u^{b}_{q}).
$$
We can also consider the translation datum $\te$ as a piecewise constant function $\te :I\to \CC$, that is we set $\te(x)=\te_{i}$ if $x\in I^{t}_{i}$. Then we call $S_{n}\te$ the Birkhoff sum over $T$ of the function $\te$. If $\ga_{q,p,n}$ is a saddle connection for $X(\pi,\la,\tau)$ in the homology class $[\ga]_{q,p,n}$ then we obviously have
\begin{equation}\label{eqholonomiaconnexsellacombinatorie}
\hol(\ga_{q,p,n})=\xi_{p}^{t}-\xi_{q}^{b}-S_{n}\te(u^{b}_{q}).
\end{equation}

\subsubsection{Reduced triples and saddle connections}\label{reducedtriplesandsaddleconnections(sss)}

Let $T:I\to I$ be an i.e.t. without connections and $q,p$ be a pair with $1<q\leq d$ and $1<p\leq d$. For any $n\in \NN$ we call $I(q,p,n)$ the open subinterval of $I$ whose endpoints are $T^{n}(u^{b}_{q})$ and $u^{t}_{p}$ (their reciprocal order does not matter).

\begin{defi}\label{deftriplaridotta}
Let $T$ and $(q,p,n)$ be respectively an i.e.t. and a triple as above. We say that $(q,p,n)$ is a \emph{reduced triple} for $T$ if for any $k\in \{0,..,n\}$ the pre-image $T^{-k}(I(q,p,n))$ of $I(q,p,n)$ does not contain in its interior any singularity $u^{t}_{p'}$ for $T$ or any singularity $u^{b}_{q'}$ for $T^{-1}$.
\end{defi}

\begin{prop}\label{propconnexsellacombinatorie}
Let $T$ be an i.e.t. without connectios defined by the data $(\pi,\la)$ (therefore $\pi$ is admissible). Let $\tau$ be any suspension datum for $\pi$ and let $X=X(\pi,\la,\tau)$ be the associated translation surface. For any triple $(q,p,n)$ reduced for $T$ there exists a saddle connection $\ga_{q,p,n}$ for $X(\pi,\la,\tau)$ which belongs to the set $\cV^{(q,p)}(\pi,\la,\tau)$ and satisfies (\ref{eqholonomiaconnexsellacombinatorie}).
\end{prop}

\begin{demo}
We can suppose without losing in generality that $T^{n}u^{b}_{q}<u^{t}_{p}$. Since $(q,p,n)$ is reduced for $T$, this is equivalent to say that $T^{k}u^{b}_{q}<T^{k-n}u^{t}_{p}$ for all $k=0,..,n$. We call $I^{(k)}$ the open interval $(T^{k}u^{b}_{q},T^{k-n}u^{t}_{p})$, so in particular $I(q,p,n)=I^{(n)}$ and all the $I^{(k)}$ have the same length. The fact that $(q,p,n)$ is reduced for $T$ means that there exists a sequence $\al(0),..,\al(n)$ of indexes in $\{1,..,d\}$ such that $I^{(k)}\subset I^{t}_{\al(k)}$ for $k=0,..,n$. In particular we observe that for any $k\in\{1,..,n\}$ we have
$$
S_{k}\te(u^{b}_{q})=S_{k}\te(T^{-n}u^{t}_{p})=
\te_{\al(0)}+..+\te_{\al(k-1)}.
$$
We set $S_{0}\te:=0$ and for $k\in\{1,..,n\}$ we introduce the simplified notation $S_{k}\te:=\te_{\al(0)}+..+\te_{\al(k-1)}$ for the sum above.

Let us fix a cartesian frame of reference on $\CC$ choosing as origin the left endpoint of the interval $I$ where $T$ acts and as positive real half-line the half-line starting from the origin and containing $I$. In particular points in $I$ are identified with complex numbers with zero imaginary part. Moreover for any $i\in\{1,..,d\}$, if $R^{t}_{i}$ and $R^{b}_{i}$ are the open rectangles introduced in paragraph \ref{veechconstruction(sss)}, we have $\overline{R^{t}_{i}}\cap \RR=I^{t}_{i}$ and $\overline{R^{b}_{i}}\cap \RR=I^{b}_{i}$

We set $v:=\xi^{t}_{p}-\xi^{b}_{q}-S_{n}\te$, then we define $\widehat{\ga}$ as the half line in the complex plane starting at $\xi^{b}_{q}$ in the direction of $v$, that is $\widehat{\ga}(t):=\xi^{b}_{q}+tv$. We observe that $\Re(v)=|I(q,p,n)|$. The path $t\mapsto \widehat{\ga}(t)$ projects to a geodesic path $t\mapsto \ga(t)$ in $X(\pi,\la,\tau)$ with $\ga(0)\in\Sigma$. A priori it is possible that $\ga$ can be extended to $\RR_{+}$, but we will show that its maximal interval of definition is $(0,1)$ with $\ga(1)\in\Sigma$, that is $\ga$ is a saddle connection.

For any $i=1,..,d$ let us call $\Pi_{i}:R^{b}_{i}\to X$ the projection from $R^{b}_{i}$ to $X$, whose image is denoted $\cR_{i}$. For $1\leq i\leq d$ the rectangles $\cR_{i}$ are disjoint each other and intersect along their boundary. Let $T>0$ be a positive real number such that $\ga$ is defined on $(0,T)$. Let $\cR_{i(0)},..,\cR_{i(k)}$ be the ordered sequence of the open rectangles defined above that $\ga(t)$ meets for $t\in(0,T)$. We observe that $\widehat{\ga}(t)$ starts in $R^{b}_{q}$, so $i(0)=q$ and for small $t$ we have $\ga(t)=\Pi_{i(0)}(\widehat{\ga}(t))$. We may also assume that $\ga(T)$ belongs to $\partial \cR_{i(k)}$, since otherwise we can extend it till to the boundary. Let $0=t_{0}<..<t_{k}<t_{k+1}=T$ be the sequence of instants such that for all $l\in\{0,..,k\}$ we have
$$
\ga(t)\in \cR_{i(l)}
\texttt{ for }
t_{l}<t<t_{l+1}
\texttt{ and }
\ga(t_{l+1})\in \partial \cR_{i(l)}.
$$

\begin{lemnonnum}
Let us suppose that the sequence $0=t_{0}<..<t_{k}<t_{k+1}=T$ satisfies $k\leq n$. Then for all $l\in\{0,..,k\}$ the following conditions hold.
\begin{description}
\item[a$(l)$]
For $t_{l}<t<t_{l+1}$ we have $\widehat{\ga}(t)+S_{l}\te\in R^{b}_{i(l)}$ and
$$
\ga(t)=\Pi_{i(l)}(\widehat{\ga}(t)+S_{l}\te).
$$
\item[b$(l)$]
$t_{l+1}<1$.
\item[c$(l)$]
$\widehat{\ga}(t_{l+1})+S_{l}\te\in I^{(l)} \subset R^{t}_{\al(l)}\cap\RR$ and there exists some $t_{l+2}>t_{l+1}$ such that $\widehat{\ga}(t)+S_{l}\te\in R^{t}_{\al(l)}$ for $t_{l+1}< t<t_{l+2}$.
%moreover for $t$ in the same interval $\ga(t)$ belongs to $\cR_{\pi(\al(l))}$ and it is given by
%$$
%\ga(t)=\rho_{\pi(\al(l))}(\widehat{\ga}(t)+S_{l+1}\te).
%$$
\end{description}
In particular we have $i(l)=\pi(\al(l-1))$ for all $l\in \{1,..,k-1\}$.
\end{lemnonnum}

\begin{demo}
The three conditions above are satisfied for $l=0$. In particular a$(0)$ is obvious, moreover $\widehat{\ga}(1)=\xi^{t}_{p}-S_{n}\te$ is not in the closure of $R^{b}_{i(0)}$, so $t_{1}<1$ and b$(0)$ also holds. Then we have $\Re(\widehat{\ga}(t_{1})-\widehat{\ga}(0))=t_{1}\Re(v)<|I^{(0)}|$, and since $I^{(0)}$ is a subinterval of $I^{b}_{i(0)}$, the condition $\ga(t_{1})\in\partial \cR_{i(0)}$ implies $\widehat{\ga}(t_{1})\in \partial R^{b}_{i(0)}\cap \RR$. In particular we have $|\widehat{\ga}(t_{1})-u^{b}_{i(0)}|<|I^{(0)}|$ and therefore $\widehat{\ga}(t_{1})$ belongs to $I^{(0)}$. Since $I^{(0)}\subset I^{t}_{\al(0)}=R^{t}_{\pi(\al(0))}\cap \RR$ then $\widehat{\ga}(t)\in R^{t}_{\al(0)}$ for $t_{1}<t<t_{2}$ and c$(0)$ follows.

Now we consider $l\in\{1,..,k\}$ and we suppose that a$(j)$,b$(j)$ and c$(j)$ hold for all $j\in\{1,..,l-1\}$. The rectangle $R^{t}_{\al(l-1)}$ is equivalent to $R^{b}_{\pi(\al(l-1))}$ via the translation by $\te_{\al(l-1)}$, so condition c$(l-1)$ implies condition a$(l)$ and we necessary have $i(l)=\pi(\al(l-1))$. Moreover $\widehat{\ga}(1)+S_{l}\te=\xi^{t}_{p}-(S_{n}\te -S_{l}\te)$, which does not belong to the closure of $R^{b}_{i(l)}$, thus $t_{l+1}<1$ and b$(l)$ follows too. Finally we have
$$
\Re(\widehat{\ga}(t_{l+1})+S_{l}\te)=
\Re(\xi^{b}_{i(0)}+S_{l}\te)+t_{l+1}\Re(v)=
T^{l}u^{b}_{i(0)}+t_{l+1}|I^{(l)}|.
$$
Condition $t_{l+1}<1$ implies that $\Re(\widehat{\ga}(t_{l+1})+S_{l}\te)\in I^{(l)}\subset I^{b}_{i(l)}$, where the inclusion of $I^{(l)}$ into $I^{b}_{i(l)}$ follows from the fact that $(q,p,n)$ is reduced. Therefore the condition $\ga(t_{l+1})\in \partial R_{i(l)}$ implies $\widehat{\ga}(t_{l+1})+S_{l}\te \in \partial R^{b}_{i(l)}\cap \RR$ and in particular $\widehat{\ga}(t_{l+1})+S_{l}\te\in I^{(l)}$. Since $I^{t}_{\al(l)}=R^{t}_{\al(l)}\cap \RR$ contains $I^{(l)}$ then c$(l+1)$ holds.

The last part of the statement is proved in the inductive proof of a$(l)$. The lemma is proved.
\end{demo}

The lemma implies that $\ga(t)$ passes through the rectangles $\cR_{i(0)},..,\cR_{i(n)}$ in $X$, where $i(0)=j$ and $i(l)=\pi(\al(l-1))$ for all $l\in\{1,..,n\}$. Then we have a sequence $0=t_{0}<..<t_{n+1}<1$ of instants such that for all $l\in\{0,..,n\}$ we have
$$
\ga(t)=\Pi_{i(l)}(\widehat{\ga}(t)+S_{l}\te)
\texttt{ and }
\ga(t_{l+1})\in\partial \cR_{i(l)}.
$$
Condition $c(n)$ says that $\widehat{\ga}(t_{n+1})+S_{n}\te$ belongs to $\partial R^{t}_{\al(n)}$, but $\ga(t_{n+1})\not \in \Sigma$. Since $\widehat{\ga}$ is a straight line in the plane with $\widehat{\ga}(1)+S_{n}\te=\xi^{t}_{p}\in\partial R^{t}_{\al(n)}$, then the instant $t_{n+2}$ in c$(n)$ equals to $1$. It follows that $\ga(t)$ can be extended for $t_{n+1}<t<1$ as the projection from $R^{t}_{\al(n)}$ to $X$ of $\widehat{\ga}(t)+S_{n}\te$. Then we evidently have $\ga(1)\in\Sigma$, thus $\ga$ is a saddle connection. It is also evident that $\hol (\ga)=v=\xi^{t}_{p}-\xi^{b}_{q}-S_{n}\te$ and that $\ga$ belongs to the set $\cV^{(q,p)}(\pi,\la,\tau)$ introduced in definition \ref{defconfigurazioneconnessionisellaveech}. The notation $\ga_{q,p,n}$ for the saddle connection $\ga$ is therefore natural. The proposition is proved.
\end{demo}

\subsection{Proof of the divergent case of theorem \ref{teoremab}}\label{prooftheoremBdivergentcase(ss)}

Let us consider a function $\varphi :[0,+\infty) \to (0,+\infty)$ such that $t\varphi(t)$ is decreasing monotone. According to remark \ref{remarkstronglydecreasingfunctions} the divergence of the integral of $\varphi$ is equivalent to the divergence of the series of the associated sequence $\varphi(n)$.

\subsubsection{Approximation of connections for i.e.t.s}\label{approximationconnectionsiet(sss)}

Let $T$ be an i.e.t. without connections. We look at triples $(q,p,n)$ with $1<q\leq d$, $1<p\leq d$ and $n\in\NN$ such that
\begin{equation}\label{eqapproxconnessionitsi}
\left|T^{n}(u_{q}^{b})-u_{p}^{t}\right|<\varphi(n).
\end{equation}

Equation (\ref{eqapproxconnessionitsi}) defines a diophantine condition for i.e.t.s inspired by the classical Khinchin condition for real numbers. Such condition has been studied in \cite{luca2}, where we prove the following result.

\begin{thm}\label{teoremaa}
Let $\varphi:[0,+\infty) \to (0,+\infty)$ be a function as above. For any admissible combinatorial datum $\pi\in S_{d}$ we have the following dichotomy:
\begin{description}
\item[a]
If $\sum_{n=1}^{+\infty}\varphi(n)<+\infty$ then equation (\ref{eqapproxconnessionitsi}) has just finitely many solutions for almost any i.e.t. $T$ in $\De_{\pi}$.
\item[b]
If $\sum_{n=1}^{+\infty}\varphi(n)=+\infty$ then for any pair $(q,p)$ with $1<q\leq d$ and $1<p\leq d$ and for almost any i.e.t. $T$ in $\De_{\pi}$ there exists infinitely many triples $(q,p,n)$ reduced for $T$ which are solution of equation (\ref{eqapproxconnessionitsi}).
\end{description}
\end{thm}

\subsubsection{Approximation of vertical saddle connections}\label{approximationconnectionssurfaces(sss)}

Let $T$ be an i.e.t. with $d$ intervals and without connections defined by the data $(\pi,\la)$. Let $\tau$ be a suspension datum for $\pi$ and $X=X(\pi,\la,\tau)$ be the associated translation structure. Recall the real vector $h\in\RR_{+}^{d}$ defined by $h_{i}=\Im(\xi^{t}_{i}-\xi^{b}_{\pi(i)})$. We consider $h$ as a piecewise constant function $h:I\to \RR_{+}$ defined by $h(x)=h_{i}$ if $x\in I^{t}_{i}$ and we denote $S_{n}h$ its Birkhoff sum.

For a pair of integers $(q,p)$ with $1<q\leq d$ and $1<p\leq d$ we consider a reduced triple $(q,p,n)$ for $T$ and the associated saddle connection $\ga_{q,p,n}$ for the surface $X=X(\pi,\la,\tau)$ given by proposition \ref{propconnexsellacombinatorie}. Since $\hol(\ga_{q,p,n})= \xi^{t}_{p}-\xi^{b}_{q}-S_{n}\te(u^{b}_{q})$ then
\begin{equation}\label{eq1approximationconnectionssurfaces(sss)}
\Re (\hol(\ga_{q,p,n}))=T^{n}u^{b}_{q}-u^{t}_{p}
\end{equation}
and we observe that $|T^{n}u^{b}_{q}-u^{t}_{p}|\leq \la_{\ast}=\sum_{k=1}^{d}\la_{k}$ for all $n\in \NN$. On the other hand
$$
\Im (\hol(\ga_{q,p,n})) =
S_{n}h(u^{b}_{q})+\Im\xi^{t}_{p}-\Im\xi^{b}_{q}
$$
and $S_{n}h(u^{b}_{q})\geq nh_{min}$ where $h_{min}:=\min_{k=1..d}h_{k}>0$. Therefore for any $\epsilon >0$ and any $n\in\NN$ big enough we have
\begin{equation}\label{eq2approximationconnectionssurfaces(sss)}
(1-\epsilon)|\ga_{q,p,n}|
\leq S_{n}h(u^{b}_{q})
\leq |\ga_{q,p,n}|.
\end{equation}

\begin{lem}\label{lemconnexsellacombinatoriequasiverticali}
Consider a function $\varphi$ as above with $\int_{0}^{\infty}\varphi(t)dt=+\infty$. Let $T$ be an i.e.t. without connections and uniquely ergodic, defined by the data $(\pi,\la)$. Let $(q,p)$ be a pair in $\{2,..,d\}^{2}$ and suppose that there exist infinitely many triples $(q,p,n)$ reduced for $T$ which are solutions of equation (\ref{eqapproxconnessionitsi}) with respect to $\varphi$. Then for any suspension datum $\tau$ in $\Th_{\pi}$ and any $\epsilon>0$ there are infinitely many saddle connections $\ga$ for the translation surface $X=X(\pi,\la,\tau)$ which belong to the set $\cV^{(q,p)}(\pi,\la,\tau)$ and satisfy
$$
|\Re(\hol(\ga))|<
\varphi\left(\frac{1-\epsilon}{1+\epsilon}\frac{|\hol(\ga)|}{\area (X)}\right).
$$
\end{lem}

\begin{demo}
According to the hypothesis we have infinitely many reduced triples $(q,p,n)$ for $T$ such that $|T^{n}u^{b}_{q}-u^{t}_{p}|<\varphi(n)$. Proposition \ref{propconnexsellacombinatorie} gives therefore infinitely many saddle connections $\ga_{q,p,n}\in\cV^{(q,p)}(\pi,\la,\tau)$ for the translation surface $X(\pi,\la,\tau)$. Since $T$ is uniquely ergodic, for $n\to \infty$ the Birkoff average $\frac{1}{n}S_{n}h(u^{b}_{q})$ converges to $\int_{I}h=\area(X)$, thus for any $\epsilon >0$ and any $\in\NN$ big enough we have
$$
(1-\epsilon)n\area(X)< S_{n}h(u^{b}_{q})<(1+\epsilon)n\area(X).
$$
Recalling that $t\varphi(t)$ is monotone decreasing, and therefore $\varphi(t)$ too, and using equations (\ref{eq1approximationconnectionssurfaces(sss)}) and (\ref{eq2approximationconnectionssurfaces(sss)}) we get the assertion. The proposition is proved.
\end{demo}

\begin{prop}\label{propdimlocaleteoremab}
Let $\varphi$ be a function as above with $\int_{0}^{\infty}\varphi(t)dt=+\infty$. Let $\pi$ be any admissible combinatorial datum in $S_{d}$ and let $(\la_{0},\tau_{0})$ be any pair of length and suspension data in $\De_{\pi}\times \Th_{\pi}$. Then there exists a neighborhood $U\subset \De_{\pi}\times \Th_{\pi}$ of $(\la_{0},\tau_{0})$ such that for any $1<q\leq d$ and $1<p\leq d$ and for almost any $(\la,\tau)\in U$ the set $\cV^{(q,p)}(\pi,\la,\tau)$ contains infinitely mant saddle connections $\ga$ for the surface $X(\pi,\la,\tau)$ such that
$$
|\Re(\hol(\ga))|<\varphi(|\hol(\ga)|).
$$
\end{prop}

\begin{demo}
Let us fix any $\epsilon>0$ and let us call $A_{0}:=\area(X(\pi,\la_{0},\tau_{0}))$. Then we consider the function
$$
\varphi_{0}(t):=
\varphi\left(\frac{1+\epsilon}{1-\epsilon}2A_{0}t\right).
$$
Our hypothesis on $\varphi$ is obviously equivalent to $\int_{0}^{\infty}\varphi_{0}(t)dt=+\infty$. Then we consider a neighborhood $\De\subset \De_{\pi}$ of $\la_{0}$ and a neighborhood $\Th\subset \Th_{\pi}$ such that for $\la\in\De$ and $\tau\in\Th$ we have
$$
A_{0}/2<\area(X(\pi,\la,\tau))<2A_{0}.
$$
The product $U:=\De\times \Th$ is obviously a neighborhood of $(\la_{0},\tau_{0})$ contained in $\De_{\pi}\times \Th_{\pi}$ and is the neighborhood were we will prove the proposition. Since $t\varphi(t)$ is decreasing monotone and therefore $\varphi(t)$ also is, for $(\la,\tau)\in U$ and for $X=X(\pi,\la,\tau)$ we have
$$
\varphi_{0}
\left(
\frac{1-\epsilon}{1+\epsilon}
\frac{t}{\area(X)}
\right)
<
\varphi(t).
$$
Theorem \ref{teoremaa} says that there exists a subset $\De'\subset \De$ of full measure such that for almost any $\la\in \De'$ and for all $1<q\leq d$ and $1<p\leq d$ there are infinitely many triples $(q,p,n)$ which are reduced for the i.e.t. $T$ defined by the data $(\pi,\la)$ and satisfy
$$
|T^{n}u^{b}_{q}-u^{t}_{p}|<\varphi_{0}(n).
$$
Moreover, according to the celebrated result of Masur and Veech (\cite{masuruno} and \cite{veech}), almost any $T$ is uniquely ergodic, therefore we can assume that for all $\la\in\De'$ the i.e.t. $T=(\pi,\la)$ is also uniquely ergodic. Lemma \ref{lemconnexsellacombinatoriequasiverticali} therefore applies and we have that for all $\la\in\De'$ and all $\tau$ in $\Th$ the set $\cV^{(q,p)}(\pi,\la,\tau)$ contains infinitely many saddle connections $\ga$ for the translation surface $X=X(\pi,\la,\tau)$ such that
$$
\Re(|\hol(\ga)|)<
\varphi_{0}
\left(
\frac{1-\epsilon}{1+\epsilon}
\frac{|\hol(\ga)|}{\area(X)}
\right)
<
\varphi(|\hol(\ga)|).
$$
We conclude observing that according to the Fubini argument $\De'\times \Th$ has full measure in $U$. The proposition is proved.
\end{demo}

\subsubsection{End of the proof}

Let $\widehat{X}_{0}=(X_{0},S_{1},..,S_{r})$ be a translation surface with frame in $\widehat{\cH}(k_{1},..,k_{r})$ such that $X_{0}$ has no vertical connection. As we explained in paragraph \ref{veechconstruction(sss)} in the background, this condition is satisfied by almost any $X_{0}\in \cH(k_{1},..,k_{r})$ and therefore for almost any $\widehat{X}_{0}$. Since $X_{0}$ has no vertical connections, there exist an admissible combinatorial datum $\pi$ in $S_{d}$, where $d-1=k_{1}+..+k_{r}$, and a pair $(\la_{0},\tau_{0})$ in $\De_{\pi}\times \Th_{\pi}$ such that $X_{0}=X(\pi,\la_{0},\tau_{0})$.

We consider the map $\cI_{\pi}:\De_{\pi}\times\Th_{\pi}\to \cH(k_{1},..,k_{r})$ introduced in paragraph \ref{veechlocalcharts(sss)} in the background. Since $\widehat{\cH}(k_{1},..,k_{r})$ is a covering space of $\cH(k_{1},..,k_{r})$, then there exist an unique continuous map
$$
\widehat{\cI}_{\pi}:\De_{\pi}\times\Th_{\pi}\to \widehat{\cH}(k_{1},..,k_{r})
$$
which lifts $\cI_{\pi}$ and such that $\widehat{\cI}_{\pi}(\la_{0},\tau_{0})=\widehat{X}_{0}$. We consider the open neighborhood $U$ of $(\la_{0},\tau_{0})$ given by proposition \ref{propdimlocaleteoremab}. Let us consider any $(\la,\tau)\in U$ such that for any $1<q\leq d$ and $1<p\leq d$, the set $\cV^{(q,p)}(\pi,\la,\tau)$ contains infinitely many saddle connections $\ga$ whose holonomy $v=\hol(\ga)$ satisfies equation (\ref{eqphiconnessionesuptrasl}). According to lemma \ref{lemconfigurazioneconnessionisellaveech}, for the corresponding $\widehat{X}=\widehat{\cI}_{\pi}(\la,\tau)$, for any pair of singular points $p_{i}$ and $p_{j}$ in $\Sigma$ and any $1\leq l\leq k_{i}$ and $1\leq m\leq k_{j}$, the set $\hol^{(p_{j},p_{i},m,l)}(\widehat{X})$ contains infinitely many solutions of equation (\ref{eqphiconnessionesuptrasl}). Corollary \ref{corveechlocalcharts}, at the level of the lift, says that $\Omega:=\widehat{\cI}_{\pi}(U)$ is an open neighborhood of $\widehat{X}_{0}$. According to proposition \ref{propdimlocaleteoremab}, the set of good $(\la,\tau)\in U$ has full measure, therefore the set of the corresponding good $\widehat{X}$ in $\Omega$ has full measure too. In other words, for a generic $\widehat{X}_{0}$, there exists a neighborhood $\Omega\subset \widehat{\cH}(k_{1},..,k_{r})$ where the statement of the diverging case of theorem \ref{teoremab} holds, therefore the statement holds almost everywhere on $\widehat{\cH}(k_{1},..,k_{r})$ and the divergent case of theorem \ref{teoremab} is proved.

\section{Asymptotic laws}

In this chapter we develop the dynamical consequences of theorem \ref{teoremab}, the dynamics we refer to being of course the Teichm\"uller flow
$$
\cF_{t}=\left(
  \begin{array}{cc}
    e^{t} & 0 \\
    0 & e^{-t} \\
  \end{array}
\right).
$$

We first recall from paragraph \ref{smoothmeasure&teichmullerflow(sss)} in the background that the group of homotheties acts on $\widehat{\cH}(k_{1},..,k_{r})$ as a subgroup of $\glduer$ by $\widehat{X}\mapsto \la\widehat{X}$ for $\la\in\RR_{+}$. The volume element on $\widehat{\cH}(k_{1},..,k_{r})$ therefore decomposes as $d\mu=\la^{n}d\la\wedge d\mu^{(1)}$, where $n=4g+2r-3$ and $d\mu^{(1)}$ is the volume element on the smooth hyperboloid $\widehat{\cH}^{(1)}(k_{1},..,k_{r})$ of area one translation surfaces with frame.

Equation (\ref{eqphiconnessionesuptrasl}), that is $|\Re(\hol (\ga))|<\varphi(|\hol(\ga)|)$, where $\ga$ is a saddle connection for a translation surface $X$, is not covariant under homotheties, indeed changing $X$ with $\la X$ it takes the form
$$
|\Re(\hol (\ga))|<\frac{\varphi(|\la\hol(\ga)|)}{\la}.
$$
Anyway, for any $\la\in\RR_{+}$, the integral $\int_{0}^{+\infty}\la^{-1}\varphi(\la t)dt$ diverges if and only if $\int_{0}^{+\infty}\varphi(t)dt$ diverges. Since the divergence or the convergence of the integral is the only discriminant condition in theorem \ref{teoremab} and as we have seen homoteties does not affect the condition, then theorem \ref{teoremab} also holds on the smooth hyperboloid $\widehat{\cH}^{(1)}(k_{1},..,k_{r})$ with respect to the measure $\mu^{(1)}$. Since the Teichm\"uller flow preserves the area of translation surfaces, it is more meaningful to study it on $\widehat{\cH}^{(1)}(k_{1},..,k_{r})$.

\subsection{Preliminary estimations.}\label{preliminaryestimations(ss)}
In this paragraph we prove some useful estimations about saddle connections which do not involve any information about the bundle they belong to, therefore in all the paragraph we will work with simple translation surfaces without any choice of a frame. Let us consider any $X$ in $\cH(k_{1},..,k_{r})$ and a saddle connection $\ga$ for $X$. As we have seen, $\ga$ is a saddle connection for $\cF_{t}X$ for any $t\in\RR$. For $t\in\RR$ let us denote $\hol(\ga,t)$ the holonomy of $\ga$ with respect to the flat structure $\cF_{t}X$. Then we set
$$
\Re(\ga,t):=\Re(\hol(\ga,t))
\texttt{  ;  }
\Im(\ga,t):=\Im(\hol(\ga,t))
\texttt{  ;  }
|\ga|_{t}:=|\hol(\ga,t)|.
$$
For simplicity we write $|\ga|$ for $|\ga|_{0}$. We also define
$$
\area(\ga,t):=|\Re(\ga,t)|\cdot|\Im(\ga,t)|
\texttt{ and }
\cot(\ga,t):=\frac{\Re(\ga,t)}{\Im(\ga,t)}.
$$
Since $\Re(\ga,t)=e^{t}\Re(\ga,0)$ and $\Im(\ga,t)=e^{-t}\Re(\ga,0)$, then $\area(\ga,t)$ is constant in $t$, on the other hand we have $\cot(\ga,t)=e^{2t}\cot(\ga,0)$. It is also easy to see that we have the relation
\begin{equation}\label{eq1sullaleggelog}
\area(\ga,t)=|\ga|_{t}^{2}\frac{|\cot(\ga,t)|}{1+\cot^{2}(\ga,t)}.
\end{equation}

\begin{lem}\label{lem1sullaleggelog}
Let $X$ be any translation surface and $\ga$ any saddle connection for $X$. If for some instant $t\geq 0$ we have $|\ga|_{t}<1$ then $t>\log |\ga|$.
\end{lem}

\begin{demo}
Since $\area(\ga,t)$ is constant and $\cot(\ga,t)=e^{2t}\cot(\ga,0)$, then equation (\ref{eq1sullaleggelog}) implies that for any $t$ we have
$$
|\ga|^{2}\frac{1}{1+\cot^{2}(\ga,0)}=
|\ga|^{2}_{t}\frac{e^{2t}}{1+\cot^{2}(\ga,t)}.
$$
We observe that $|\cot(\ga,t)|>|\cot(\ga,0)|$, since $t\geq 0$. Then $|\ga|_{t}<1$ implies that $e^{2t}>|\ga|^{2}$ and the lemma follows.
\end{demo}

Let us fix any $X\in\cH(k_{1},..,k_{r})$ and any saddle connection $\ga$ for $X$. Let us denote $\tau:=\tau(X,\ga)\in \RR$ the instant defined by
$$
|\ga|_{\tau}:=\min \{|\ga|_{t};t\in\RR\}.
$$
Since the geometry of a translation surface is locally euclidian, the length $|\ga|_{t}$ is minimal when $\Re(\ga,t)=\Im(\ga,t)$, that is $\cot(\ga,t)=1$ and we have
\begin{equation}\label{eq2sullaleggelog}
\tau(X,\ga)=-\frac{1}{2}\log(\cot(\ga,0)).
\end{equation}
Since by definition $|\cot(\ga,\tau)|=1$, recalling that $\area(\ga,t)$ is constant, equation (\ref{eq1sullaleggelog}) implies
\begin{equation}\label{eq3sullaleggelog}
|\ga|^{2}_{\tau}=2\area(\ga,0)
\texttt{ and }
|\ga|^{2}_{t}>2\area(\ga,0)
\texttt{ for }
t\not =\tau.
\end{equation}

\begin{lem}\label{lem2sullaleggelog}
Let $\epsilon$ be any positive real number. For almost any $X$ in $\cH(k_{1},..,k_{r})$ and for all but finitely many saddle connection $\ga$ for $X$, the instant $\tau(X,\ga)$ satisfies:
$$
\tau(X,\ga)\leq (1+\epsilon)\log |\ga|.
$$
\end{lem}

\begin{demo}
Since the function $t\mapsto t^{-(1+2\epsilon)}$ has convergent tail, the convergent case of theorem \ref{teoremab} implies that for almost any $X\in\cH(k_{1},..,k_{r})$ and for all but finitely many saddle connections $\ga$ for $X$ we have $|\Re(\hol(\ga))|\geq |\ga|^{-(1+2\epsilon)}$, which implies
$$
\cot(\ga,0)\geq
\frac{1}{|\Im(\ga,0)||\ga|^{1+2\epsilon}}\geq
|\ga|^{-2(1+\epsilon)}.
$$
Then equation (\ref{eq2sullaleggelog}) implies that for such $X$ and for such $\ga$ the condition in the statement holds. The lemma is proved.
\end{demo}

\subsection{Proof of theorem \ref{teoremac}.}\label{prooftheoremc(ss)}

In this section we prove theorem \ref{teoremac}. Before entering into details we state the following useful lemma, whose proof is a trivial exercise in calculus.

\begin{lem}\label{lemprooftheoremc}
Let $\varphi:[1,+\infty)\to (0,+\infty)$ be a function such that $t\varphi(t)$ is decreasing monotone. The function $\widehat{\varphi}:[0,+\infty)\to (0,+\infty)$ defined by $\widehat{\varphi}(s):=e^{s}\varphi(e^{s})$ is decreasing monotone.

On the other hand, for any decreasing monotone function $\psi:[0,+\infty)\to (0,+\infty)$ there exists an unique function $\varphi:[1,+\infty)\to (0,+\infty)$ such that $t\varphi(t)$ is decreasing monotone and such that $\psi=\widehat{\varphi}$. The function $\varphi$ is given by $\varphi(t)=\psi(\log t)/t$ and we have
$$
\int_{0}^{\infty}\psi(t)dt=\int^{\infty}_{1}\varphi(s)ds.
$$
\end{lem}

\subsubsection{Convergent case}

Here we consider the case $\int_{0}^{\infty}\psi(t)dt<+\infty$.
%The function $X\mapsto\sys(X)$ does not need any choice of frame of horizontal separatirces, indeed it equals to the minimum value of $\sys^{(p_{i},l,p_{j},m)}(\widehat{X})$ for all pair of points $p_{i},p_{j}$ in $\Sigma$ and all pair of indexes $1\leq i\leq k_{i}$ and $1\leq m\leq k_{j}$ and it just depends from $X$. Therefore we can forget about the information on the frame of horizontal separatirces and look at orbits $\cF_{t}X$ of elements $X\in\cH^{(1)}(k_{1},..,k_{r})$.
Our aim is to prove that for generic $X$ we have
$$
\lim_{t\to\infty}\frac{\sys(\cF_{t}X)}{\sqrt{\psi(t)}}=\infty.
$$
Our strategy is first to prove that for generic $X$ we have
\begin{equation}\label{eq1proofoftheoremteoremac}
\liminf_{t\to\infty}\frac{\sys(\cF_{t}X)}{\sqrt{\psi(t)}}>1.
\end{equation}
Then, once (\ref{eq1proofoftheoremteoremac}) is proved, we observe that for any positive constant $C>1$ the function $C^{2}\psi$ still has finite integral, thus we can change $\psi$ with $C^{2}\psi$ and (\ref{eq1proofoftheoremteoremac}) still holds for $C^{2}\psi$, which is equivalent to say that $\liminf_{t\to\infty}\sys(\cF_{t}X)/\sqrt{\psi(t)}>C$. We get that
$\lim_{t\to\infty}\sys(\cF_{t}X)/\sqrt{\psi(t)}$ exists and is equal to $+\infty$.

We pass to the proof of equation (\ref{eq1proofoftheoremteoremac}). Let us suppose that there exists a positive measure subset $\cS$ of $\cH^{(1)}(k_{1},..,k_{r})$ such that for any $X\in \cS$ we have
$$
\liminf_{t\to\infty}
\frac{\sys(\cF_{t}X)}
{\sqrt{\psi(t)}}
\leq 1.
$$
Let us fix any $X\in\cS$. There exists a sequence of instants $t_{1}<...<t_{n}<...$ with $t_{n}\to +\infty$ and such that $\sys(\cF_{t_{n}}X)\leq \sqrt{\psi(t_{n})}$. Let $\ga_{n}$ be a sequence of saddle connections for $X$ such that for every $n\in \NN$ we have $\sys(\cF_{t_{n}}X)=|\ga_{n}|_{t_{n}}$. For any such $\ga_{n}$ we obviously have $|\ga_{n}|^{2}_{t_{n}}\leq \psi(t_{n})$. Since $\psi$ has convergent tail, in particular is goes to zero at infinity and it follows that $|\ga_{n}|_{t_{n}}\to 0$ as $n\to\infty$. We observe that $|\ga_{n}|_{t}$ as $t\to+\infty$ for any fixed $n$. Since $t_{n}\to +\infty$, then any $\ga_{n}$ occurs in the sequence $(\ga_{n})$ finitely many times. Moreover on any translation surface $X$ there are just finitely many saddle connections with length smaller than some fixed bound, thus it follows that $|\ga_{n}|_{0}\to\infty$ as $n\to\infty$ (i.e. the length at $t=0$). We also observe that this implies that $\cot(\ga_{n},0)\to 0$ as $n\to\infty$. Indeed if we could find a positive constant $c$ (depending on $X$) and a subsequence $(n(k))_{k\in\NN}$ with $|\cot(\ga_{n(k)},0)|>c$, then equation (\ref{eq1sullaleggelog}) would imply $\area(\ga_{n(k)},t_{n(k)})\to +\infty$ as $k\to\infty$, but equation (\ref{eq3sullaleggelog}) implies that this is in contradiction with our assumption on $\ga_{n}$.

Let us fix any $\epsilon>0$. Since $\cot(\ga_{n},0)\to 0$ then for any $n$ big enough we have
$$
|\Im(\ga_{n},0)|>(1+\epsilon)\frac{\area(\ga_{n})}{|\ga_{n}|_{0}}
$$
and since by equation (\ref{eq3sullaleggelog}) we have $\area(\ga_{n},0)\leq (1/2)|\ga_{n}|^{2}_{t_{n}}$ then it follows that
$$
|\Re(\hol(\ga_{n}))|=|\Re(\ga_{n},0)|<
\frac{1+\epsilon}{2}\frac{|\ga_{n}|^{2}_{t_{n}}}{|\ga_{n}|_{0}}<
\frac{1+\epsilon}{2}\frac{\psi(t_{n})}{|\ga_{n}|_{0}}.
$$
Finally, recalling lemma \ref{lem1sullaleggelog}, we have $t_{n}>\log(|\ga_{n}|_{0})$ and it follows that
$$
|\Re(\hol(\ga_{n}))|<
\frac{1+\epsilon}{2}\frac{\psi(\log(|\ga_{n}|_{0}))}{|\ga_{n}|_{0}}.
$$
We set $\varphi(t):=\psi(\log t)/t$ and according to lemma \ref{lemprooftheoremc} we have $\int_{1}^{+\infty}\varphi(t)dt<+\infty$. On the other hand, for a set $\cS$ of positive measure in $\cH(k_{1},..,k_{r})$ and for all $X$ in $\cS$ we can find infinitely many saddle connection $\ga_{n}$ for $X$ such that $|\Re(\hol(\ga_{n}))|<\varphi(|\hol(\ga_{n})|)$, which is absurd, since is in contradiction with the convergent part of theorem \ref{teoremab}. Equation (\ref{eq1proofoftheoremteoremac}) is proved an therefore the convergent part of theorem \ref{teoremac} too.

\subsubsection{Divergent case.}

Now we consider the case $\int_{0}^{\infty}\psi(t)dt=+\infty$. We recall that for $\widehat{X}\in \widehat{\cH}^{(1)}(k_{1},..,k_{r})$ the bundles $\cC^{(p_{j},p_{i},m,l)}(\widehat{X})$ are invariant for the Teichm\"uller flow $\cF_{t}$ (lemma \ref{lemconfigurazioneconnessionisella}). Let us fix any pair of points $p_{j},p_{i}\in\Sigma$ and any pair of indexes $m,l$ with $m\in\{1,..,k_{j}\}$ and $l\in\{1,..,k_{i}\}$. We prove that for almost any $\widehat{X}\in \widehat{\cH}^{(1)}(k_{1},..,k_{r})$ we have
$$
\liminf_{t\to\infty}
\frac{\sys^{(p_{j},p_{i},m,l)}(\cF_{t}\widehat{X})}
{\sqrt{\psi(t)}}
=0.
$$
Similarly to the convergent case, our strategy is to prove that for generic $\widehat{X}$ we have
\begin{equation}\label{eq2proofoftheoremteoremac}
\liminf_{t\to\infty}
\frac{\sys^{(p_{j},p_{i},m,l)}(\cF_{t}\widehat{X})}
{\sqrt{\psi(t)}}
<1.
\end{equation}
Once (\ref{eq2proofoftheoremteoremac}) is proved, we observe that for any positive constant $\epsilon>0$ the function $\epsilon^{2}\psi$ still has divergent integral. Therefore we can change $\psi$ with $\epsilon^{2}\psi$ and (\ref{eq2proofoftheoremteoremac}) still holds for $\epsilon^{2}\psi$, which is equivalent to
$$
\liminf_{t\to\infty}
\frac{\sys^{(p_{j},p_{i},m,l)}(\cF_{t}\widehat{X})}
{\sqrt{\psi(t)}}
<\epsilon.
$$
Since the argument works for all $\epsilon$ we get $\liminf_{t\to\infty}\sys^{(p_{j},p_{i},m,l)}(\cF_{t}\widehat{X})\sqrt{\psi(t)}=0$.

We pass to the proof of (\ref{eq2proofoftheoremteoremac}). Let us fix any $\epsilon>0$ and let us consider the function
$$
\varphi(s):=
\frac{\psi((1+\epsilon)\log (s))}{s}
$$
Lemma \ref{lemprooftheoremc} (modulo a trivial change of variable) implies $\int_{1}^{+\infty}\varphi(t)dt=+\infty$, therefore the divergent part of theorem \ref{teoremab} applies. For almost any $\widehat{X}\in\widehat{\cH}^{(1)}(k_{1},..,k_{r})$ the bundle $\cC^{(p_{j},p_{i},m,l)}(\widehat{X})$ contains infinitely many saddle connections $\ga_{n}$ such that $|\Re(\hol(\ga_{n}))|<\varphi(|\hol(\ga_{n})|)$. It follows that we have
$$
\area(\ga_{n})<
|\Im(\ga_{n},0)|\varphi(|\ga_{n}|)<
|\ga_{n}|\varphi(|\ga_{n}|),
$$
where $|\ga_{n}|=|\hol(\ga_{n})|$ is the initial length $|\ga_{n}|_{0}$. For any such $\ga_{n}$ let us consider the instant $\tau(n)$ such that $|\ga_{n}|_{\tau(n)}=\min \{|\ga_{n}|_{t};t\in\RR\}$. According to equation (\ref{eq3sullaleggelog}), the length $|\ga|_{\tau(n)}$ of the saddle connection $\ga_{n}$ with respect to $\cF_{\tau(n)}X$ satisfies
$$
|\ga_{n}|^{2}_{\tau(n)}<(1/2)|\ga_{n}|\varphi(|\ga_{n}|).
$$
Since we are considering a generic $\widehat{X}$ we may assume that the underlying translation surface $X$ belongs to the full measure subset of those elements for which lemma \ref{lem2sullaleggelog} applies. For the length of $\ga_{n}$ at $t=0$ we have $|\ga_{n}|\geq\exp \left(\frac{\tau(n)}{1+\epsilon}\right)$
and since the function $t\varphi(t)$ is decreasing monotone we get
$$
|\ga_{n}|^{2}_{\tau(n)}<
\frac{1}{2}
\exp\left(\frac{\tau(n)}{1+\epsilon}\right)
\varphi\left(
\exp\left(\frac{\tau(n)}{1+\epsilon}\right)
\right)
=\psi(\tau(n)).
$$
Since the bundles of saddle connections are preserved by the action of $\cF_{t}$ on $\widehat{\cH}^{(1)}(k_{1},..,k_{r})$ then for any $\tau(n)$ the saddle connection $\ga_{n}$ belongs to the bundle $\cC^{(p_{j},p_{i},m,l)}(\cF_{\tau(n)}\widehat{X})$, therefore
$$
\sys^{(p_{j},p_{i},m,l)}(\cF_{\tau(n)}\widehat{X})\leq |\ga_{n}|_{\tau_{n}}<\sqrt{\psi(\tau(n))}.
$$
Finally we observe that $|\Re(\ga_{n},0)|\to 0$ and $|\ga_{n}|_{0}\to \infty$ as $n\to \infty$, therefore $\cot(\ga_{n},0)\to 0$ and equation (\ref{eq2sullaleggelog}) implies that $\tau_{n}\to\infty$. Equation (\ref{eq2proofoftheoremteoremac}) therefore follows and the divergent part of theorem \ref{teoremac} is proved.

\section{Punctured tori}

This section is devoted to the proof of theorem \ref{teoremad}. Let us consider a function $\varphi:[0,+\infty)\to(0,\infty)$ such that $t\varphi(t)$ is decreasing monotone and $\int_{0}^{+\infty}\varphi(t)dt=+\infty$. Let $X$ be any flat torus with $r$ punctures $p_{1},..,p_{r}$. We want to prove that for almost any $\te\in [0,2\pi)$ there exist at least a pair of points $(p_{j},p_{i})$ such that $\hol^{(p_{j},p_{i})}(X_{\te})$ contains infinitely many solutions of equation (\ref{eqphiconnessionesuptrasl}).

\subsection{Veech's construction for punctured flat tori}\label{veechconstructionflattori(ss)}

\subsubsection{First return to an horizontal section}\label{firstreturnhorizontalsection(sss)}

Let $X$ be a flat torus with $r$ punctures and let $\Sigma=\{p_{1},..,p_{r}\}$ be the set of the punctures. We suppose that $X$ does not have vertical saddle connections. Let $I$ be an horizontal segment in $X$ whose endpoints lie on two trajectories of the vertical field $\partial_{y}$, respectively ending and starting at the same point $p_{\ast}\in\Sigma$. It is a well known fact that the first return $T$ of the vertical flow to $I$ is an irrational rotation. Let us consider the length coordinate $x$ on $I$, which give an identification of $I$ with the interval $(0,|I|)\subset\RR$. There exists some $\de\in(0,|I|)$ such that $T$ acts by
$$
T(x)=x+\de
\texttt{ if }
0<x<|I|-\de
$$
$$
T(x)=x-(|I|-\de)
\texttt{ if }
|I|-\de<x<|I|.
$$
In particular $T$ is an i.e.t. and $(0,|I|-\de)$ and $(|I|-\de,|I|)$ are the two maximal intervals where $T$ acts as a translation. Anyway, since we want to keep track of all points in $\Sigma$, we consider a partition of $I$ into $d=r+1$ intervals $I^{t}_{1},..,I^{t}_{d}$, any $I^{t}_{i}$ being delimited by two endpoints whose positive trajectory under the vertical flow ends in some point of $\Sigma$. Therefore $T$ is described by a pair of data $(\pi,\la)$, where $\la_{i}=|I^{t}_{i}|$ for $i\in\{1,..,d\}$ and $\pi\in S_{d}$ is a \emph{rotational} combinatorial datum, that is it satisfies
$$
\pi(i)-i=const \mod d.
$$
With this notation $T$ has a real discontinuity at the point on $I$ whose coordinate is $u^{t}_{\pi^{-1}(1)}$, and is continuous at all $u^{t}_{i}$ with $i\not =\pi^{-1}(1)$. Similarly $T^{-1}$ is continuous at all points $u^{b}_{j}$ with $j\not =\pi^{1}(1)$ and has a real discontinuity just at $u^{b}_{\pi(1)}$. We define the normalized length datum $\widehat{\la}:=\frac{\la}{\la_{\ast}}$, where $\la_{\ast}=\sum_{i=1}^{d}\la_{i}$, and we consider the i.e.t. $\widehat{T}$ defined by the data $(\pi,\widehat{\la})$. According to the discussion above, $\widehat{T}$ is a rotation and it has two maximal intervals where it acts as a translation, anyway we look at all points $\widehat{u}^{t}_{p}$ and $\widehat{u}^{b}_{q}$ with $1<q\leq d$ and $1<p\leq d$. In terms of the data $(\pi,\la)$ they are given by
$$
\widehat{u}^{t}_{p}=\sum_{j<p}\widehat{\la}_{j}
\texttt{ and }
\widehat{u}^{b}_{q}=\sum_{\pi(j)<\pi(q)}\widehat{\la}_{j}.
$$
In particular the \emph{rotation number} of $\widehat{T}$ is $\al=\al(\la):=\widehat{u}^{b}_{\pi(1)}$ and $\widehat{T}$ acts by
$$
\widehat{T}(x)=x+\al
\texttt{ if }
0<x<1-\al
$$
$$
T(x)=x-(1-\al)
\texttt{ if }
1-\al<x<1.
$$

We can reconstruct the original torus $X$ from $T$ with the Veech construction, that is there exists some $\tau\in\De_{\pi}$ such that $X=X(\pi,\la,\tau)$. In general, any translation surface obtained with the Veech construction from a rotational combinatorial datum $\pi\in S_{d}$ is a flat torus with $d-1$ marked points.

\subsubsection{Section for the rotated vertical flow}\label{sectionrotatedverticalflow(sss)}

Let $\pi\in S_{r+1}$ be a rotational combinatorial datum and for $(\la,\tau)\in \De_{\pi}\times \Th_{\pi}$ let $X=X(\pi,\la,\tau)$ be the associated flat torus with $r$ marked points. For $\te$ in $\RR$ we define the vectors in $\RR^{d}$
$$
\la(\te):=\cos\te \la -\sin\te\tau
\texttt{ and }
\tau(\te):=\sin\te\la+\cos\te\tau.
$$
There exists an open interval $J=J(\pi,\la,\tau)\subset \RR$ containing $0$ and depending on $\pi,\la, \tau$ such that $\la(\te)\in\De_{\pi}$ and $\tau(\te)\in\Th_{\pi}$ for $\te\in J$. For $\te$ in the same interval the rotated torus $X_{\te}$ can be therefore represented with the Veech construction as $X_{\te}=X(\pi,\la(\te),\tau(\te))$. For $\te\in R$ we define $\la_{\ast}(\te)=\sum_{i=1}^{d}\la_{i}(\te)$ and $\widehat{\la}(\te):=\frac{\la(\te)}{\la_{\ast}(\te)}$ and we consider the i.e.t. $\widehat{T}_{\te}$ whose data are $(\pi,\widehat{\la}(\te))$. It is easy to check that the map $\te\mapsto \widehat{\la}(\te)$ is the parametrization of a segment in the standard simplex $\De^{(1)}:=\{\la\in\RR^{d}_{+};\la_{\ast}=1\}$, that is we can write
$$
\widehat{\la}(\te)=\widehat{\la}+s(\te)v,
$$
where $s:J\to \RR$ is a smooth change of variable and $v$ is a vector in $\RR^{d}$. Their explicit expression is
$$
s(\te)=\frac{\tan\te}{\la_{\ast}-\tau_{\ast}\tan\te}
\texttt{ and }
v=\frac{\tau_{\ast}}{\la_{\ast}}\la-\tau,
$$
where $\tau_{\ast}=\sum_{i=1}^{d}\tau_{i}$. It follows that for $\te\in J(\pi,\la,\tau)$ and for $1<q\leq d$ and $1<p\leq d$ the marked points $\widehat{u}^{t}_{p}(\te)$ and $\widehat{u}^{b}_{q}(\te)$ in $(0,1)$ are given by
$$
\widehat{u}^{t}_{p}(\te)=
\widehat{u}^{t}_{p}(0)+s(\te)\sum_{k<p}v_{k}
\texttt{ and }
\widehat{u}^{b}_{q}(\te)=
\widehat{u}^{b}_{q}(0)+s(\te)\sum_{\pi(k)<\pi(q)}v_{k}.
$$
According to the discussion above the only real singularities for $\widehat{T}_{\te}$ and $\widehat{T}^{-1}_{\te}$ are respectively $\widehat{u}^{t}_{\pi^{-1}(1)}(\te)$ and $\widehat{u}^{b}_{\pi(1)}(\te)$. In particular, observing that $\la_{\ast}\sum_{\pi(j)\leq \pi(1)}v_{j}=\area(X)$, the rotation number $\al(\te)$ of $\widehat{T}_{\te}$ is given by
$$
\al(\te)=\widehat{u}^{b}_{\pi(1)}(\te)=
\al(0)+s(\te)\frac{\area(X)}{\la_{\ast}}.
$$
Let us fix any triple of data $(\pi,\la,\tau)$. We resume the discussion of these two paragraphs in the following lemma

\begin{lem}\label{lemsingolaritatsitororotato}
The map $\te\mapsto \al(\te)$, which assigns to any $\te\in J(\pi,\la,\tau)$ the rotation number $\al(\te)$ of $\widehat{T}_{\te}$, is a smooth change of variable. Moreover, for any $1<i\leq d$ there exists two pairs of real numbers $(A^{t}_{i},B^{t}_{i})$ and $(A^{b}_{i},B^{b}_{i})$, depending only from $(\pi,\la,\tau)$, such that
$$
\widehat{u}^{t}_{i}(\te)=A^{t}_{i}+\al(\te)B^{t}_{i}
\texttt{ and }
\widehat{u}^{b}_{i}(\te)=A^{b}_{i}+\al(\te)B^{b}_{i}.
$$
\end{lem}

\subsection{Proof of theorem \ref{teoremad}}\label{prooftheoremD(ss)}

Let us fix a rotational $\pi\in S_{d}$ and a pair of length and suspension data $\la\in\De_{\pi}$ and $\tau\in\Th_{\pi}$. Let $J=J(\pi,\la,\tau)$ be the interval such that $(\la_{\te},\tau_{\te})\in\De_{\pi}\times\Th_{\pi}$ for all $\te\in J$. For any such $\te$ let $T_{\te}$ be the i.e.t. whose data are $(\pi,\la_{\te})$ and for a pair $(q,p)$ with $1<q\leq d$ and $1<p\leq d$ let $u^{b}_{q}(\te)$ and $u^{t}_{p}(\te)$ be the corresponding singularities. We look at triples $(q,p,n)$ satisfying
\begin{equation}\label{eqtriplenonridottegenereuno}
|T^{n}_{\te}u^{b}_{q}(\te)-u^{t}_{p}(\te)|<\varphi(n).
\end{equation}
We call $\cP_{\te}$ the subset of $\{2,..,d\}^{2}$ of those pairs $(q,p)$ as above such that there exist infinitely many triples $(q,p,n)$ which are reduced for $T_{\te}$ and satisfy (\ref{eqtriplenonridottegenereuno}).

\begin{lem}\label{lemtriplenonridottegenereuno}
For almost any $\te\in J$ and for any pair $(q,p)$ with $1<q\leq d$ and $1<p\leq d$ there exist infinitely many triples $(q,p,n)$ which satisfy equation (\ref{eqtriplenonridottegenereuno}), anyway these triples are not necessarily reduced for $T_{\te}$.
\end{lem}

\begin{demo}
Let us fix any pair $(q,p)$ with $1<q\leq d$ and $1<p\leq d$. Since $\widehat{T}_{\te}$ is a rotation on $(0,1)$ with rotation number $\al(\te)$, for any $x\in(0,1)$ and any $\in\NN$ we have $\widehat{T}^{n}_{\te}(x)=\{x+n\al(\te)\}$, where $\{\cdot\}$ denotes the fractionary part. Then we have
$$
|T^{n}_{\te}u^{b}_{q}(\te)-u^{t}_{p}(\te)|=
\la_{\ast}(\te)
|\widehat{T}^{n}_{\te}\widehat{u}^{b}_{q}(\te)-
\widehat{u}^{t}_{p}(\te)|=
\la_{\ast}(\te)|\{\widehat{u}^{b}_{q}(\te)+n\al(\te)\}-
\widehat{u}^{t}_{p}(\te)|.
$$
Since $\varphi(n)\to 0$ as $n\to\infty$, for any $\te\in J$, in order to have infinite solutions $n\in\NN$ of equation (\ref{eqtriplenonridottegenereuno}), it is enough to find infinitely many solutions of
$$
\{\widehat{u}^{b}_{q}(\te)-\widehat{u}^{t}_{p}(\te)+n\al(\te)\}<
\frac{\varphi(n)}{\la_{\ast}(\te)}.
$$
Let us consider the decomposition $J=\bigsqcup_{k\in\ZZ}J_{k}$, where $2^{k}\leq \la_{\ast}(\te)<2^{k+1}$ for $\te\in J_{k}$. Let us fix $k\in\ZZ$ and let us consider the function $\varphi_{k}:=2^{-(k+1)}\varphi$. Using the formula in lemma \ref{lemsingolaritatsitororotato} for $\widehat{u}^{b}_{q}(\te)$ and $\widehat{u}^{t}_{p}(\te)$ it is evident that it is sufficient to have, for any $\te\in J_{k}$, infinitely many solutions of
$$
\{(B^{b}_{q}-B^{t}_{p}+n)\al(\te)-(A^{b}_{q}-A^{t}_{p})\}<
\varphi_{k}(n).
$$
Since $\int_{0}^{\infty}\varphi_{k}(t)dt=\infty$, theorem \ref{teoremae} in paragraph \ref{arithmeticresult(ss)} implies that there exist infinitely many $n\in\NN$ satisfying the inequality above for almost any $\al$, moreover according to lemma \ref{lemsingolaritatsitororotato}, the change of variable $\al=\al(\te)$ is smooth and in particular it preserves sets of measure zero. It follows that we have infinitely many solutions for almost any $\te\in J_{k}$ too. Since the result holds for any $k$ and the union is countable, the required condition holds for almost any $\te\in J$. The lemma is proved.
\end{demo}

We need the following lemma, which holds in general for any i.e.t. and not just for rotational ones.

\begin{lem}\label{lemtriplenonridottegenerale}
For any i.e.t. $T$ without connections there exists $\epsilon>0$ such that the following is true. If the triple $(q,p,N)$ is not reduced for $T$ and satisfies
$$
|T^{N}u^{b}_{q}-u^{t}_{p}|< \epsilon
$$
then there exist $q'$ and $p'$ with $1<q'\leq d$ and $1<p'\leq d$ and positive integers $n$ and $m$ with $n,m<N$ such that the triples $(q,p',n)$ and $(q',p,m)$ are reduced for $T$ and satisfy
$$
|T^{n}u^{b}_{q}-u^{t}_{p'}|<|T^{N}u^{b}_{q}-u^{t}_{p}|
\texttt{ and }
|T^{m}u^{b}_{q'}-u^{t}_{p}|<|T^{N}u^{b}_{q}-u^{t}_{p}|.
$$
\end{lem}

\begin{demo}
We call $u(1)<..<u(2d-2)$ the elements in the set $\{u^{t}_{i};1<i\leq d\}\sqcup\{u^{b}_{j};1<j\leq d\}$, displayed in increasing order. Then we put $2\epsilon:=\min_{l=2,..,2d-2}|u(i)-u(i-1)|$, which is positive since $T$ has not connections. In particular our choice implies that all the intervals $I^{t}_{i}$ and $I^{b}_{j}$ have length at least $2\epsilon$. Let $(q,p,N)$ be a triple as in the statement and let $I(q,p,N)$ be the interval whose endpoints are $T^{N}_{\te}(u^{b}_{q})$ and $u^{t}_{p}$. Since $(q,p,N)$ is not reduced for $T$ there exists some $k\in \{0,..,n\}$ and $l\in\{2,..,d\}$ such that $T^{-k}(I(q,p,N))$ contains in its interior either $u^{t}_{l}$ or $u^{b}_{l}$, moreover it follows by our choice of $\epsilon$ that in fact $k \geq 1$. We consider the smallest $k$ such that the last condition holds. By minimality we have that $T^{-i}$ is a translation on $I(q,p,N)$ for any $i=0,..,k$, that is $|T^{N-k}u^{b}_{q}-T^{-k}u^{t}_{p}|=|T^{N}u^{b}_{q}-u^{t}_{p}|<\epsilon$. Without any loss in generality we suppose that we have $u^{t}_{l}\in T^{-k}(I(q,p,n))$. We can also suppose that we have
$$
T^{N-k}u^{b}_{q}<u^{t}_{l}<T^{-k}u^{t}_{p}.
$$
We first look at the inequality on the right. By our choice of $\epsilon$, no other singularity of $T$ or $T^{-1}$ is contained in the interval $(u^{t}_{l},T^{-k}u^{t}_{p})$ and by minimality of $k$ this also holds for the iterates $T^{i}(u^{t}_{l},T^{-k}u^{t}_{p})$ for $i=0,..,k$. The interval $(u^{t}_{l},T^{-k}u^{t}_{p})$ is mapped by $T$ onto $(u^{b}_{l},T^{-(k-1)}u^{t}_{p})$. Then we apply $T$ again $k-1$ times and we get that the triple $(l,p,k-1)$ is reduced and satisfies$|T^{k-1}u^{b}_{l}-u^{t}_{p}|<|T^{N}u^{b}_{q}-u^{t}_{p}|$.

Now we look at the inequality on the left in the condition above. It implies that the triple $(q,l,N-k)$ satisfies $|T^{N-k}u^{b}_{q}-u^{t}_{l}|<|T^{N}u^{b}_{q}-u^{t}_{p}|$. If $(q,l,N-k)$ is not reduced for $T$, we call $j(q,p,N)$ the number of singularities, both for $T$ and $T^{-1}$, which are contained in the orbit
$$
I(q,p,N),T^{-1}(I(q,p,N)),..,T^{-(N-k)}(I(q,p,N)).
$$
We observe that $j(q,l,N-k)<j(q,p,N)$, therefore we can start a descending induction procedure until we get a reduced triple. The lemma is proved.
\end{demo}

\begin{cor}\label{cortriplenonridottegenereuno}
For almost any $\te\in J$ the set $\cP_{\te}$ contains at least $2r-1$ elements.
\end{cor}

\begin{demo}
Let us fix any $\te$ in the full measure subset of $J$ given by lemma \ref{lemtriplenonridottegenereuno}. Let $(q,p)$ be any pair in $\{2,..d\}^{2}$ and according to lemma \ref{lemtriplenonridottegenereuno} consider a family of infinitely many triples $(q,p,N_{k})$ satisfying (\ref{eqtriplenonridottegenereuno}). According to lemma \ref{lemtriplenonridottegenerale}, either there exists a subsequence $N_{i}$ of $N_{k}$ such that the corresponding triples $(q,p,N_{i})$ are reduced for $T_{\te}$, that is $(q,p)\in\cP_{\te}$, or there are two indexes $q',p'$ with $1<q'\leq d$ and $1<p'\leq d$ with $(q,p')\in\cP_{\te}$ and $(q',p)\in\cP_{\te}$. We display the elements of $\{2,..,d\}$ in a $r\times r$ matrix. The argument above implies that any row and any column contain at least an element of $\cP_{\te}$, that is the latter has at least $2r-1$ elements.
\end{demo}

\subsubsection{End of the proof}

Let $X_{0}$ be any flat torus with $r$ marked points and let $X$ be an element in $\soduer X_{0}$ without saddle connections.

According to the discussion in paragraph \ref{firstreturnhorizontalsection(sss)}, let us chose a rotational $\pi\in S_{d}$ and length-suspension data $(\la,\tau)\in\De_{\pi}\times \Th_{\pi}$ in order to have $X=X(\pi,\la,\tau)$. Then consider the open interval $J=J(\pi,\la,\tau)$ and for $\te$ in the full measure subset of $J(\pi,\la,\tau)$ given by corollary \ref{cortriplenonridottegenereuno} consider the rotated torus $X_{\te}=X(\pi,\la_{\te},\tau_{\te})$ and the associated i.e.t. $T_{\te}$ whose combinatorial and length data are $(\pi,\la_{\te})$. Observe that we obviously have $\area(X_{\te})=\area(X)$. Let us fix any $\epsilon>0$ and consider the function
$$
\varphi_{\epsilon}(t):=
\varphi\left(
\frac{1+\epsilon}{1-\epsilon}\area(X)t
\right).
$$
Since $\varphi_{\epsilon}$ differs form $\varphi$ just for a linear change of variable, then it satisfies the same properties, that is $t\varphi_{\epsilon}(t)$ is decreasing monotone and $\int_{0}^{\infty}\varphi_{\epsilon}(t)dt=\infty$.

For any pair $(q,p)$ in $\cP_{\te}$ there are infinitely many $n\in\NN$ such that the corresponding triples $(q,p,n)$ are reduced for $T_{\te}$ and satisfy equation (\ref{eqtriplenonridottegenereuno}). Lemma \ref{lemconnexsellacombinatoriequasiverticali} in paragraph \ref{approximationconnectionssurfaces(sss)} implies that the set $\cV^{(q,p)}(\pi,\la_{\te},\tau_{\te})$ contains infinitely many saddle connections $\ga$ for $X_{\te}$ such that
$$
|\Re(\hol(\ga))|<
\varphi_{\epsilon}\left(
\frac{1-\epsilon}{1+\epsilon}\frac{|\hol(\ga)|}{\area (X)}
\right)
=\varphi(|\hol(\ga)|).
$$
According to corollary \ref{cortriplenonridottegenereuno} the set $\cP_{\te}$ contains at least $2r-1$ pairs $(q,p)$. Lemma \ref{lemconfigurazioneconnessionisellaveech} in paragraph \ref{relationveech(sss)} implies that to any such $(q,p)$ it corresponds an unique pair $(p_{j},p_{i})\in\Sigma^{2}$ such that the set $\cV^{(q,p)}(\pi,\la_{\te},\tau_{\te})$ coincides with $\cC^{(p_{j},p_{i})}(X_{\te})$. The set $\hol^{(p_{j},p_{i})}(X_{\te})$ therefore contains infinitely many solutions of equation (\ref{eqphiconnessionesuptrasl}) and we have at least $2r-1$ different sets like $\hol^{(p_{j},p_{i})}(X_{\te})$.

We argue that the set of $X$ without saddle connections has full measure in $\soduer X_{0}$. Since for any such $X$ exists an interval $J_{X}\subset \RR$ containing zero such that $X_{\te}$ satisfies the condition in theorem \ref{teoremad} for almost any $\te\in J_{X}$, then the condition holds almost everywhere on $\soduer X_{0}$. Theorem \ref{teoremad} is proved.

\subsection{Appendix: an arithmetic result}\label{arithmeticresult(ss)}

We consider a function $\varphi:[0,+\infty)\to(0,\infty)$ such that $t\varphi(t)$ is decreasing monotone and any pair $(x,y)$ of real numbers. For any $\al\in(0,1)$ we look for solutions $n\in\NN$ of the following equation
\begin{equation}\label{eqteoremae}
\{(n+x)\al-y\}<\varphi (n).
\end{equation}

\begin{thm}\label{teoremae}
Let $\varphi$ be a function as above.
\begin{itemize}
\item
If $\sum_{n=1}^{\infty}\varphi(n)<+\infty$, then for almost every $\al \in \RR_{+}$ there exist finitely many solutions $n\in \NN$ of equation (\ref{eqteoremae}).
\item
If $\sum_{n=1}^{\infty}\varphi(n)=\infty $ then for almost any $\al$ equation (\ref{eqteoremae}) has infinitely many solutions $n\in\NN$.
\end{itemize}
\end{thm}

The convergent part of the theorem it is a straightforward application of the easy half of the \emph{Borel-Cantelli lemma}, anyway we give here a short proof. We fix any pair of positive integers $m$ and $n$ and we introduce the function $f_{n,m}$ from $[m,m+1)$ to $[0,1)$, defined by $f_{n,m}(\al):=\{(n+x)\al-y\}$, which is piecewise linear and has the same slope on all its branches. The branches of $f_{n,m}$ are defined on sub-intervals of $[m,m+1)$, if we exclude the leftmost and the rightmost of these intervals all the others have the same length (equal to $(n+x)^{-1}$) and restricted to them $f_{n,m}$ is surjective onto $[0,1)$. It follows than for any $\epsilon>0$ and any $n$ big enough, for any subinterval $I$ of $[0,1)$ we have $\leb(f^{-1}_{n,m}(I))\leq (1+\epsilon)\leb(I)$. In particular for any $m$ in $\NN$ and for any $n$ big enough we have
$$
\leb(f^{-1}_{n,m}(0,\varphi(n)))\leq (1+\epsilon)\varphi (n).
$$
Since $\varphi$ has convergent series it follows that $\sum_{n=1}^{\infty}\leb(f^{-1}_{n,m}(0,\varphi(n)))< +\infty$ and accordingly to the Borel-Cantelli lemma, equation (\ref{eqteoremae}) has just finitely many solutions for almost any $\al$ in $[m,m+1)$. The same argument works on any other interval $[m,m+1)$ and the proof of the convergent part is complete. The proof of the divergent part of theorem \ref{teoremae} is more complicated and its proof takes the remaining of this section.

\subsubsection{Notation}

For any positive real number $\al$ we introduce the symbol $\widehat{\al}$ to denote the vector $(1,\al)\in \RR_{+}^{2}$, that is the vector with unitary horizontal component and slope equal to $\al$. For any pair $v,w\in \RR^{2}_{+}$ of linear independent vectors in the first quadrant we denote $[v,w]$ the matrix in $\glduer$ whose first and second columns are respectively $v$ and $w$. We define a total ordering on $\PP\RR^{2}_{+}$ writing for any such a pair $v,w$
$$
v\prec w
\textrm{  if  }
\det [v,w]>0
\texttt{ and }
v\preceq w
\textrm{  if  }
\det [v,w]\geq 0.
$$

\subsubsection{Classical continued fraction}\label{classicalcontinuedfraction(sss)}

In this paragraph we recall some classical facts about the \emph{continued fraction algorithm}, which gives good rational approximations $p/q$ of a real number $\al$. We will follow a geometrical interpretation. We are interested to real numbers admitting an infinite sequence of approximations, therefore we just consider irrational $\al$. In this case, for any pair $p,q$ of positive integers, denoting $r=(q,p)\in\NN^{2}$, we have
$$
\det[r,\widehat{\al}]\not =0.
$$

Let us set $r_{-2}:=(1,0)\in \NN^{2}$ and $r_{-1}:=(0,1)\in \NN^{2}$. For any $\al\in \RR_{+}$ we define a vector $r_{0}=r_{0}(\al)\in \NN^{2}$ by
$$
r_{0}(\al):=a_{0}(\al)r_{-1}+r_{-2}
$$
where $a_{0}(\al)\in \NN$ is such that $a_{0}r_{-1}+r_{-2}\preceq(1,\al)\prec (a_{0}+1)r_{-1}+r_{-2}$. Let us write the vector $r_{0}(\al)$ as $r_{0}(\al)=(q_{0}(\al),p_{0}(\al))$, with $q_{0}(\al),p_{0}(\al)\in\NN$. Letting $\al$ vary we introduce the family of integer vectors $Q_{0}:=\left\{r_{0}(\al);\al\in \RR_{+}\setminus \QQ\right\}$ and the partition $\cQ_{0}:=\left\{I(r_{0});r_{0}\in Q_{0}\right\}$ whose elements are the intervals $I(r_{0})$ with constant value for the function $\al\mapsto r_{0}(\al)$.

Now let us suppose that for any $\al\in \RR_{+}\setminus \QQ$ and for all $i<n$ we have defined all the vectors $r_{i}(\al)=(q_{i}(\al),p_{i}(\al))$, with $q_{i}(\al),p_{i}(\al)\in\NN$, the associated families of integer vectors $Q_{i}:=\left\{r_{i}(\al);\al\in \RR_{+}\setminus\QQ\right\}$ and the partitions $\cQ_{i}:=\left\{I(r_{i});r_{i}\in Q_{i}\right\}$ whose elements are the intervals $I(r_{i})\subset\RR_{+}$ with specified value for the function $\al\mapsto r_{i}(\al)$. Then we define by induction
$$
r_{n}(\al):=a_{n}(\al)r_{n-1}(\al)+r_{n-2}(\al)
$$
where $a_{n}(\al)$ is the positive integer defined according to the following condition
$$
a_{n}r_{n-1}+r_{n-2}\prec \widehat{\al} \prec (a_{n}+1)r_{n-1}+r_{n-2}
\textrm{  if  }n \textrm{  is even  }
$$
$$
(a_{n}+1)r_{n-1}+r_{n-2}\prec \widehat{\al} \prec a_{n}r_{n-1}+r_{n-2}
\textrm{  if  }n \textrm{  is odd  }.
$$
We write $r_{n}(\al)=(q_{n}(\al),p_{n}(\al))$ with $q_{n}(\al),p_{n}(\al)\in\NN$ and letting $\al$ vary we define the family of integer vectors $Q_{n}:=\left\{r_{n}(\al);\al\in \RR_{+}\setminus \QQ\right\}$ and the sigma-algebra $\cQ_{n}$, whose atoms are the intervals $I(r_{n})\subset\RR_{+}$ with specified value for the function $\al\mapsto r_{n}(\al)$. For any fixed $\al\in \RR_{+}\setminus \QQ$ the vectors $r_{n}(\al)$ are a sequence of approximations of $\widehat{\al}$ in $\PP\RR^{2}_{+}$. The corresponding sequence $p_{n}/q_{n}$ of rational numbers is the sequence of approximations of $\al$ with respect to the continued fraction algorithm. We recall some well known properties for these approximations that we use in the following, proofs can be found in \cite{kin}.
\begin{enumerate}
\item
For any $\al$ and any $n$ we have $r_{2n}\prec \widehat{\al}\prec r_{2n+1}$. The even approximations $r_{2n}$ give a monotone increasing sequence and the odd ones a monotone decreasing sequence. Moreover, if $n$ is even and if $r_{n}$ is any element of $Q_{n}$, there exist a decreasing sequence $(r_{n+1}(a))_{a\in\NN}$ of elements of $Q_{n+1}$ and an increasing sequence $(r_{n+2}(a))_{a\in\NN}$ of elements of $Q_{n+2}$ such that respectively $r_{n+1}(a)\to r_{n}$ and $r_{n+2}(a)\to r_{n}$ in $\PP\RR^{2}_{+}$ as $a\to \infty$. The analogous property holds for $n$ odd.
\item
For any $n\in \NN$ let us introduce the vector $r_{n}'=r_{n}'(\al):=r_{n}(\al)+r_{n-1}(\al)$ and let us write it as $r'_{n}=(q'_{n},p'_{n})$ with $q'_{n},p'_{n}\in\NN$. We have$|\det[r_{n},r'_{n}]|=1$, that is $\{r_{n},r'_{n}\}$ is a basis of $\ZZ^{2}$.
\item
For any $n$ and any $r_{n}$ we have $|I(r_{n})|=(q_{n}q'_{n})^{-1}$ and $q^{2}_{n+2}>2q^{2}_{n}$. We also have the uniform estimation $|I(r'_{n})|<|I(r_{n})|<3|I(r'_{n})|$ for two consecutive atoms $I(r_{n})$ and $I(r'_{n})$ of $\cQ_{n}$.
\item
There exists a constant $\ga>1$ such that for almost any $\al$ and any $n$ big enough the denominators $q_{n}=q_{n}(\al)$ of the approximations of $\al$ satisfy $q_{n}<e^{\ga n}$.
\end{enumerate}

The family of sigma-algebras $\{\cQ_{n}\}_{n\in\NN}$ satisfies $\cQ_{n}\subset\cQ_{n+1}$ for any $n\in\NN$, that is it defines a \emph{monotone increasing filtration}, moreover $\sup_{n}\cQ_{n}$ is the Borel sigma-algebra. It follows that any interval $J$ contained in $\RR_{+}$, with its endpoints included or not, admits a decomposition $J=\bigsqcup_{n=0}^{\infty}J_{n}$, where $J_{0}$ is the maximal subset of $J$ measurable with respect to $\cQ_{0}$ and $J_{n}$ is defined iteratively for $n\geq1$ as the maximal subset of $J\setminus\bigsqcup_{i=0}^{n-1}J_{i}$ which is measurable with respect to $\cQ_{n}$. To any interval $J$ we associate an integer index $i(J)$ defined as the minimum of those $n\in\NN$ such that $J_{n}\not = \emptyset$. For any $n>i(J)$ we introduce the sets $G_{n}=\bigsqcup_{l=i(J)}^{n}J_{l}$ and $\cG_{n}:=J\setminus G_{n}$. Applying the properties (1) and (3) of the continued fraction algorithm listed above it is easy to get the following lemma, whose proof is left to the reader.

\begin{lem}\label{lemclassicalcontinuedfraction}
There exist a positive $\la$ with $0<\la<1$ such that for any interval $J\subset\RR_{+}$ and for any $n> i(J)$ the decomposition $J=G_{n}\sqcup\cG_{n}$ satisfies the following properties
\begin{enumerate}
\item
$G_{n}$ is measurable with respect to $\cQ_{n}$ and $|\cG_{n}|\leq \la^{n-i(J)}|J|$.
\item
If $I(r_{n})$ is any atom of $\cQ_{n}$ contained in $G_{n}$, we have $|I(r_{n})|<\la^{n-i(J)}|J|$.
\item
If $J'$ is any non-empty connected component of $\cG_{n}$, we have $i(J')\leq n+2$.
\end{enumerate}
\end{lem}

\subsubsection{Twisted continued fraction.}\label{twistedcontinuedfraction(sss)}

In this paragraph we develop some useful machinery to solve equation \ref{eqteoremae}, our exposure is based on a construction given in \cite{ydue} at pages 105,106.

We fix a pair of real numbers $(x,y)$ and we look for approximations of $\al$ of the form $\frac{j+y}{k+x}$ with $j,k\in\ZZ^{2}$. We use the same vectorial notation of paragraph \ref{classicalcontinuedfraction(sss)}, thus we introduce the vector $v=(x,y)\in\RR^{2}$ and for any $r=(k,j)\in\ZZ^{2}$ we consider the vector $s=r+v=(k+x,j+y)$, which are elements of the coset $\ZZ^{2}+v$ of $\ZZ^{2}$ in $\RR^{2}$. We assume that $\al$ is irrational and moreover that for any $s\in\ZZ^{2}+v$ we have
$$
\det[\widehat{\al},s]\not =0.
$$
We denote $\cA$ the set of those $\al$ in $\RR_{+}\setminus\QQ$ such that the condition above is satisfied. It is evidently a full measure subset of $\RR_{+}$. Let us fix $\al\in\cA$ and let $(r_{k}(\al))_{k\in\NN}$ be the sequence of its classical approximations introduced in paragraph \ref{classicalcontinuedfraction(sss)}. For any $k\in \NN$ let us define
$$
\Lambda(r_{k}):=\{sr_{k}+tr'_{k};s,t\in(0,1]\},
$$
that is the parallelogram in $\RR^{2}$ spanned by the pair of vectors $(r_{k},r_{k}')$, which is a fundamental domain for the action of $\slduez$ on $\RR^{2}$. Then we define $v(r_{k})$ as the only element in $(\ZZ^{2}+v)\cap \Lambda(r_{k})$.

We observe that for $k$ even the $k$-th approximation $r_{k}(\al)=(q_{k}(\al),p_{k}(\al))$ satisfies $r_{k}\prec \widehat{\al}\prec r'_{k}$. In order to define the sequence of the \emph{twisted approximations} $s_{n}(\al)$ of $\al$ we introduce a parameter $N\in\NN^{\ast}$ and for $n\in\NN$ we consider the $2N(n-1)$-th approximation $r_{2N(n-1)}=r_{2N(n-1)}(\al)$ of $\al$ with respect to the continued fraction algorithm of paragraph \ref{classicalcontinuedfraction(sss)}. We observe that we have $r_{2N(n-1)}\prec v(r_{2N(n-1)})\prec r'_{2N(n-1)}$ and we define $s_{n}$ according to the following two cases:

\begin{itemize}
\item
If $r_{2N(n-1)}\prec \widehat{\al} \prec v(r_{2N(n-1)})$ we define $\nu_{n}=\nu_{n}(\al)$ as the minimum of those $\nu$ in $\NN$ such that $v(r_{2N(n-1)})+\nu r_{2N(n-1)}\prec \widehat{\al}$. In this case we always have $\nu_{n}(\al)\geq1$. Then we define the $n$-th twisted approximation as
\begin{equation}\label{eq1twistedcontinuedfraction}
s_{n}:=v(r_{2N(n-1)})+\nu_{n}r_{2N(n-1)}.
\end{equation}
In this case we also define the vector $s'_{n}=s'_{n}(\al)$ by $s'_{n}:=s_{n}-r_{2N(n-1)}$.
\item
If $v(r_{2N(n-1)})\prec \widehat{\al}\prec r'_{2N(n-1)}$ we define $\nu_{n}=\nu_{n}(\al)$ as the maximum of those $\nu$ in $\NN$ such that $v(r_{2N(n-1)})+\nu r'_{2N(n-1)}\prec \widehat{\al}$. In this case we may also have $\nu_{n}(\al)=0$. Then we define the $n$-th twisted approximation as
\begin{equation}\label{eq1(bis)twistedcontinuedfraction}
s_{n}:=v(r_{2N(n-1)})+\nu_{n}r'_{2N(n-1)}.
\end{equation}
In this case we also define the vector $s'_{n}=s'_{n}(\al)$ by $s'_{n}:=s_{n}+r'_{2N(n-1)}$.
\end{itemize}

Observe that by definition we have $s_{n}(\al)\prec \widehat{\al}\prec s'_{n}(\al)$ for any $n\in \NN$. We write $s_{n}(\al)=(k_{n}(\al)+x,j_{n}(\al)+y)$ with  $k_{n}(\al),j_{n}(\al)\in\ZZ$ and similarly $s'_{n}=(k'_{n}+x,j'_{n}+y)$, with $k'_{n},j'_{n}\in\ZZ$. If $s_{n}(\al)$ is defined by equation (\ref{eq1twistedcontinuedfraction}) we have $\det[s_{n},s'_{n}]=\det[r_{2N(n-1)},v(r_{2N(n-1)})]$. On the other hand, if $s_{n}(\al)$ is defined by equation (\ref{eq1(bis)twistedcontinuedfraction}) we have $\det[s_{n},s'_{n}]=\det[v(r_{2N(n-1)}),r'_{2N(n-1)}]$. In both cases, recalling that $\det[r_{2N(n-1)},r'_{2N(n-1)}]=1$ and that $v(r_{2N(n-1)})$ belongs to the fundamental domain spanned by these two vectors, we have
\begin{equation}\label{eq2twistedcontinuedfraction}
0<\det[s_{n},s'_{n}]\leq1.
\end{equation}

Letting $\al$ vary we define the family of vectors $P_{n}:=\{s_{n}(\al);\al\in\cA\}$, which is a subset of $(\ZZ^{2}+v)\cap\RR_{+}^{2}$, and the sigma-algebra $\cP_{n}$, whose atoms are the intervals $I(s_{n})\subset\RR_{+}$ with specified value for the function $\al\mapsto s_{n}(\al)$. We observe that for any $n\in \NN$ the $\si$-algebra $\cP_{n}$ is a refinement of $\cQ_{2N(n-1)}$, anyway the family $\{\cP_{n}\}_{n\in\NN^{\ast}}$ it is not a filtration.

\begin{lem}\label{lem2twistedcontinuedfraction}
For any $n\in \NN$ and for any atom $I(s_{n})$ of $\cP_{n}$ we have $i(I(s_{n}))\leq 2N(n-1)+3$.
\end{lem}
\begin{demo}
We treat separately the two cases where $s_{n}$ is defined according to equation (\ref{eq1twistedcontinuedfraction}) or equation (\ref{eq1(bis)twistedcontinuedfraction}). If $r_{2N(n-1)}\prec s_{n}\prec v(r_{2N(n-1)})$ we observe that there exists $a\in\NN^{\ast}$ such that
$$
s_{n}\prec r'_{2N(n-1)}+ar_{2N(n-1)}\prec s'_{n}.
$$
For the real numbers $\al\in I(r_{2N(n-1)})$ with $a_{2N(n-1)+1}(\al)=a+1$ we have $r'_{2N(n-1)}+ar_{2N(n-1)}=r'_{2N(n-1)-1}+(a+1)r_{2N(n-1)}=r_{2N(n-1)+1}(\al)$, that is $I(s_{n})$ contains an element of $\cQ_{2N(n-1)+1}$. It follows from the properties of the continued fraction that there exist a sequence $r_{2N(n-1)+2}$ of elements of $\cQ_{2N(n-1)+2}$ which accumulate to $r_{2N(n-1)+1}$. The corresponding interval $I(r_{2N(n-1)}+2)$ are therefore eventually contained in $I(s_{n})$ and thus $i(I(s_{n}))=2N(n-1)+2$.

If $v(r_{2N(n-1)})\preceq s_{n}\prec r'_{2N(n-1)}$ we observe that there exists $a\in\NN^{\ast}$ such that
$$
s_{n}\prec r_{2N(n-1)}+ar'_{2N(n-1)}\prec s'_{n}.
$$
For the real numbers $\al\in I(r_{2N(n-1)})$ with $a_{2N(n-1)+1}(\al)=1$ and $a_{2N(n-1)+2}(\al)=a$ we have respectively $r'_{2N(n-1)}=r_{2N(n-1)+1}(\al)$ and $r_{2N(n-1)}+ar'_{2N(n-1)}=r_{2N(n-1)+2}(\al)$, that is $I(s_{n})$ contains an element of $\cQ_{2N(n-1)+2}$. Arguing as in the previous case we get $i(I(s_{n}))=2N(n-1)+3$. the lemma is proved.
\end{demo}

\subsubsection{A sufficient condition}\label{sufficientcondition(cap5)(sss)}

For any $n\in\NN^{\ast}$ and any $s_{n}\in\cP_{n}$ we define $\Upsilon[s_{n}]$ as the bijective affine function from $I(s_{n})$ to $[0,1)$, or in formula:
$$
\Upsilon[s_{n}](\al):=
\frac{k'_{n}+x}{\det[s_{n},s'_{n}]}\det[s_{n},\hat{\al}].
$$
For any $\al\in I(s_{n})$ we have by definition $s_{n}\prec \widehat{\al}\prec s'_{n}$, thus equation (\ref{eq2twistedcontinuedfraction}) implies $0<\det[s_{n},\widehat{\al}]<\det[s_{n},s'_{n}]\leq1$. Since $\det[s_{n},\widehat{\al}]=(k_{n}(\al)+x)\al-(j_{n}(\al)+y)$ we get $\det[s_{n},\widehat{\al}]=\{(k_{n}+x)\al-y\}$. It follows that in order to have infinitely many solution of equation (\ref{eqteoremae}) it is enough to have infinitely many solutions $n$ of
$$
\Upsilon[s_{n}(\al)](\al)<\frac{k'_{n}+x}{\det[s_{n},s'_{n}]}\varphi(k_{n}).
$$
As a consequence of lemma \ref{lem2twistedcontinuedfraction} we have $k_{n}+x=\Re(s_{n})<\Re(r_{2N(n-1)+3})=q_{2N(n-1)+3}$ and similarly for $k'_{n}$. Since there exists a constant $\ga>1$ such that for almost any $\al$ and any $n$ big enough we have $q_{n}<e^{\ga n}$, then we can find a constant $\ga'>1$ such that for almost any $\al$ and any $n$ big enough we have
$$
k_{n}+x<e^{\ga'n}
\texttt{ and }
k'_{n}+x<e^{\ga'n}.
$$
We define a positive sequence setting $\psi_{n}:=e^{\ga'n}\varphi(e^{\ga'n})$. Since $t\varphi(t)$ is decreasing monotone, then $\psi_{n}$ is decreasing monotone and for almost any $\al$ and any $n$ big enough we have $\psi_{n}<(k'_{n}(\al)+x)\varphi(k_{n}(\al))$. Moreover according to remark \ref{remarkstronglydecreasingfunctions} condition $\int_{0}^{+\infty}\varphi(t)dt=+\infty$ is equivalent to $\sum_{n=1}^{+\infty}\psi_{n}=\infty$. Finally, according to equation (\ref{eq2twistedcontinuedfraction}) we have $0<\det[s_{n},s'_{n}]<1$, therefore in order to have infinitely many solution of equation (\ref{eqteoremae}) for almost any $\al$ it is enough to have infinitely many solutions $n$ of
\begin{equation}\label{eqsufficientcondition(cap5)}
\Upsilon[s_{n}(\al)](\al)<\psi_{n}.
\end{equation}
It follows that in order to prove the divergent part of theorem \ref{teoremae} it is enough to prove the following proposition

\begin{prop}\label{propsufficientcondition(cap5)}
If $\psi:\NN\to \RR_{+}$ is a decreasing monotone sequence such that $\sum_{n\in\NN}\psi_{n}=+\infty$ then for almost every $\al \in \RR_{+}$ there are infinitely many solutions $n\in\NN$ of (\ref{eqsufficientcondition(cap5)}).
\end{prop}

\subsubsection{End of the proof}\label{endproof(cap5)(sss)}

In this paragraph we prove proposition \ref{propsufficientcondition(cap5)}. Since any function $\Upsilon[s_{n}]$ has image in $(0,1)$ we can assume that there exists some positive $\epsilon\in(0,1)$ such that we definitively have $\psi_{n}<1-\epsilon$, otherwise the statement is trivially true. Then we make a choice of the parameter $N$ in the definition of the twisted approximations requiring that $N\geq 3$ and $\la^{5-2N}<\epsilon$, where $\la$ is the constant in lemma \ref{lemclassicalcontinuedfraction}.\\

For any $n\in\NN^{\ast}$ and any $s_{n}\in\cP_{n}$ we define $J(s_{n})$ as the subinterval of $I(s_{n})$ of those $\al$ with $\Upsilon[s_{n}](\al)\geq\psi_{n}$. We evidently have $|J(s_{n})|= (1-\psi_{n})|I(s_{n})|$.

\begin{lem}\label{lemendoftheproof(cap5)}
For any $n\in\NN^{\ast}$ and any $s_{n}\in\cP_{n}$ we have $i(J(s_{n}))<2Nn$.
\end{lem}

\begin{demo}
If for some $s_{n}$ in some $\cP_{n}$ we have $i(J(s_{n}))\geq 2Nn$, then $J(s_{n})$ cannot contain elements $r_{2Nn-2}$ of $Q_{2Nn-2}$, otherwise a sequence of intervals $I(r_{2Nn-1})$ would be contained in $J(s_{n})$. It follows that there exists some atom of $\cQ_{2Nn-2}$ which contains $J(s_{n})$ as subinterval. On the other hand $i(I(s_{n}))\leq 2N(n-1)+3$, according to lemma \ref{lem2twistedcontinuedfraction} and $2Nn-2-i(I(s_{n}))\geq 2N-5\geq 1$. Then we can apply the decomposition in lemma  \ref{lemclassicalcontinuedfraction} and we write $I(s_{n})=G_{2Nn-2}\sqcup\cG_{2Nn-2}$. Let $I(r_{2Nn-2})$ be the atom of $\cQ_{2Nn-2}$ which contains $J(s_{n})$. Either $I(r_{2Nn-2})$ is contained in $G_{2Nn-2}$ or $J(s_{n})$ is a subset of $\cG_{2Nn-2}$. In both cases lemma \ref{lemclassicalcontinuedfraction} implies that  $|J(s_{n})|<\la^{2N-5}|I(s_{n})|$. Finally we recall that have $|J(s_{n})|>\epsilon|I(s_{n})|$, because $\psi_{n}<1-\epsilon$, thus we get to an absurd since $\la^{2N-5}<\epsilon$. The lemma is proved.
\end{demo}

We set $\cC_{n}:=\bigcup_{s_{n}\in\cP_{n}}J(s_{n})$ and $\cC:=\bigcup_{m\in\NN}\bigcap_{n\geq M}\cC_{n}$. The statement in proposition \ref{propsufficientcondition(cap5)} corresponds to $\leb (\cC)=0$ and in order to show it it is sufficient to prove that for any $M\in \NN$ we have $\leb(\bigcap_{n\geq M}\cC_{n})=0$. Thus we fix any $M\in\NN$ and we re-define $\cC:=\bigcap_{n\geq M}\cC_{n}$. To show that $\cC$ has zero measure we define a family of sets $\widehat{\cC}_{n}$ with $\widehat{\cC}_{n}\supset \widehat{\cC}_{n+1}$ and such that $\cC\subset \widehat{\cC}_{n}$ for any $n\in\NN$ and $\leb(\widehat{\cC}_{n})\to 0$ as $n\to\infty$. The definition is given with the following induction procedure.

\begin{itemize}
\item
The first element of the family corresponds to $n=M$ and we set $\widehat{\cC}_{M}:=\cC_{M}$. Any connected components of $\widehat{\cC}_{M}$ is a sub-interval $J(s_{M})$ of some atom $I(s_{M})$ of $\cP_{M}$ and according to lemma \ref{lemendoftheproof(cap5)} we have $i(J(s_{M}))< 2NM$. Therefore lemma \ref{lemclassicalcontinuedfraction} implies that we can decompose it as
$$
J(s_{M})=G_{2NM}\sqcup\cG_{2NM},
$$
where the subset $G_{2NM}$ is measurable with respect to the sigma-algebra $\cQ_{2NM}$, and therefore with respect to $\cP_{M+1}$, and where the remaining part satisfies $|\cG_{2NM}|<\la|J(s_{M})|$ ($\la$ is the constant appearing in lemma \ref{lemclassicalcontinuedfraction}). Moreover, if $J'$ is any non-empty connected component of $\cG_{2NM}$ we have $i(J')\leq2NM+2$.
\item
We fix $n>M$ and suppose that the sets $\widehat{\cC}_{M},..,\widehat{\cC}_{n-1}$ are defined. We also assume that for any $k\in \{M,..,n-1\}$ any connected component $J$ of $\widehat{\cC}_{k}$ satisfies $i(J)< 2Nk$ and we observe that the assumption is satisfied for $k=M$.

Let $J$ be any connected component of $\widehat{\cC}_{n-1}$. Since $i(J)<2N(n-1)$, using lemma \ref{lemclassicalcontinuedfraction} we decompose it as
$$
J=G_{2N(n-1)}\sqcup\cG_{2N(n-1)},
$$
where the subset $G_{2N(n-1)}$ is measurable with respect to the sigma-algebra $\cQ_{2N(n-1)}$, and therefore with respect to $\cP_{n}$, and where the remaining part satisfies $|\cG_{2N(n-1)}|<\la|J(s_{n})|$ and $i(J')\leq2N(n-1)+2$ for any non-empty connected component $J'$ of $\cG_{2N(n-1)}$. We define
$$
G_{2N(n-1)}\cap\widehat{\cC}_{n}:=G_{2N(n-1)}\cap\cC_{n},
$$
that is we consider the union of the subintervals $J(s_{n})$ of those $I(s_{n})$ contained in $G_{2N(n-1)}$. According to lemma \ref{lemendoftheproof(cap5)} any new connected component $J(s_{n})$ of $\widehat{\cC}_{n}$ arising in this way satisfies $i(J(s_{n}))< 2Nn$. Then we complete the definition setting
$$
\cG_{2N(n-1)}\cap\widehat{\cC}_{n}:=\cG_{2N(n-1)},
$$
that is the rest $\cG_{2N(n-1)}$ passes unchanged to $\widehat{\cC}_{n}$. According to lemma \ref{lemclassicalcontinuedfraction}, for these connected component of $\widehat{\cC}_{n}$, i.e. for the non empty connected components $J'$ of $\cG_{2N(n-1)}$ we have $i(J')\leq 2N(n-1)+2<2Nn$.
\end{itemize}

For any $n\geq M$ the sets $\widehat{\cC}_{n}$ is therefore defined and if $J$ is any of its connected component we have $i(J)<2Nn$.

\begin{lem}\label{lem1thesufficientconditionhastotalmeasure(cap5)}
For any $n\geq M$ and any connected component $J$ of $\widehat{\cC}_{n}$ we have
$$
|J\cap \widehat{\cC}_{n+1}|<(1-\de\psi_{n+1})|J|,
$$
where $\de=1-\la$ and $\la\in (0,1)$ is the constant in lemma \ref{lemclassicalcontinuedfraction}.
\end{lem}

\begin{demo}
Let $J$ be a connected component of $\widehat{\cC}_{n}$. Since $i(J)<2Nn$, we apply lemma \ref{lemclassicalcontinuedfraction} and we decompose it as $J=G_{2Nn}\sqcup \cG_{2Nn}$. According to the inductive definition of $\widehat{\cC}_{n+1}$ we have
$$
J\cap\widehat{\cC}_{n+1}=(G_{2Nn}\cap\cC_{n+1})\sqcup \cG_{2Nn}.
$$
Since $G_{2Nn}$ is union of atoms $I(s_{n+1})$ of $\cP_{n+1}$ and for any such atom we have $|I(s_{n+1})\cap\cC_{n+1}|=|J(s_{n+1})|=(1-\psi_{n+1})|I(s_{n+1})|$, we get $|G_{2Nn}\cap\cC_{n+1}|=(1-\psi_{n+1})|G_{2Nn}|$. On the other hand, lemma \ref{lemclassicalcontinuedfraction} also implies $|\cG_{2Nn}|<\la|J|$, therefore summing the two contributions we have
$$
|J\cap \widehat{\cC}_{n+1}|=
(1-\psi_{n+1})|G_{2Nn}|+|\cG_{2Nn}|=
\left((1-\psi_{n+1})\frac{|G_{2Nn}|}{|J|}+\frac{|\cG_{2Nn}|}{|J|}\right)|J|=
$$
$$
\left(1-(1-\frac{|\cG_{2Nn}|}{|J|})\psi_{n+1}\right)|J|<
(1-(1-\la)\psi_{n+1})|J|.
$$
\end{demo}

Now let us consider any connected component $J$ of $\widehat{\cC}_{M}$. Lemma \ref{lem1thesufficientconditionhastotalmeasure(cap5)} implies that for any $n>M$ we have
$$
|J\cap \widehat{\cC}_{n}|<
\prod_{k=M+1}^{n}(1-\de\psi_{k})|J|.
$$
For any $M\geq 1$, the condition $\sum_{n\in\NN}\psi_{n}=+\infty$ on the sequence $\psi_{n}$ is equivalent to $\lim_{n\to\infty}\prod_{k=M+1}^{n}(1-\de\psi_{k})=0$. It follows that $|J\cap \widehat{\cC}_{n}|\to 0$ as $n\to\infty$ for any connected component $J$ of $\widehat{\cC}_{M}$. Therefore $|J\cap \cC|= 0$ for any $J$. Since $\widehat{\cC}_{M}$ is countable union of its connected component we get $|\cC|=0$. Proposition \ref{propsufficientcondition(cap5)} is proved and therefore the divergent part of theorem \ref{teoremae} too.

\end{document}